\documentclass[final,1p,times]{elsarticle}
\usepackage{algorithm}
\usepackage{algorithmic}
\usepackage[algo2e]{algorithm2e}
\usepackage{amssymb}
\usepackage{amsthm}
\usepackage{tikz}
\usepackage{amscd}
\usepackage{amsmath}
\usepackage{amsfonts}
\usepackage{booktabs}
\usepackage{graphicx}
\usepackage{color}
\usepackage{framed}
\usepackage{comment}
\usepackage{chngcntr}
\usepackage{mathrsfs}
\usepackage{float}
\usepackage{subcaption}
\usepackage{hyperref}
\usepackage[pagewise,mathlines]{lineno}
\usepackage{geometry}

\numberwithin{equation}{section}
\usepackage{bm}
\usepackage{bbding}

\newtheorem{theorem}{Theorem}[section]

\newtheorem{lemma}[theorem]{Lemma}

\makeatletter

\newcommand{\Rmnum}[1]{\expandafter\@slowromancap\romannumeral #1@}
\long\def\tbl#1#2{\begin{center}\caption{#1}#2\end{center}}
\makeatother
\newcommand{\colrule}{\midrule}
\newcommand{\botrule}{\bottomrule}

\biboptions{sort&compress}
\journal{***}

\begin{document}
\begin{sloppypar}
\begin{frontmatter}
\title{Evolution of weights on a connected finite graph}
\author[author1]{Jicheng Ma}
\ead{2019202433@ruc.edu.cn}
\author[author3]{Yunyan Yang{\footnote{corresponding author}}}
\ead{yunyanyang@ruc.edu.cn}
\address{$^1$School of Mathematics, Renmin University of China, Beijing, 100872, China}
\begin{abstract}
On a connected finite graph, we propose an evolution of weights including
Ollivier's Ricci flow as a special case. During the evolution process, on each edge, the speed of change of weight is exactly the difference between the Wasserstein distance related to two probability measures and certain graph distance. Here the probability measure may be chosen as
an $\alpha$-lazy one-step random walk, an $\alpha$-lazy two-step random walk, or a general probability measure. Based on the ODE theory, we show that the initial value problem has a unique global solution.

A discrete version of the above evolution is applied to the problem of community detection. Our algorithm is based on such a discrete evolution,
where probability measures are chosen as $\alpha$-lazy one-step random walk and $\alpha$-lazy two-step random walk respectively. Note that the latter measure has not been used in previous works. Moreover, only one surgery needs to be performed after the last iteration, which makes our algorithm much simpler than earlier ones based on Lin-Lu-Yau's Ricci curvature.
\end{abstract}
\begin{keyword}
weighted graph; Ollivier's Ricci flow; Lin-Lu-Yau's Ricci flow; community detection
\\
\MSC[2020] 05C21; 05C85; 35R02; 68Q06
\end{keyword}
\end{frontmatter}
\section{Introduction}

In 1982, Ricci flow was  first introduced by Hamilton \cite{Hamilton} to deform a metric on a Riemannian manifold $(M,g)$ according to the differential equation involving a one-parameter family of metrics $g(t)$ on a time interval, namely
$$
\partial_t g(t)=-2{\rm Ric}(g(t)),
$$
 where ${\rm Ric}(g(t))$ is the Ricci curvature of the metric $g(t)$. Generally speaking,
Ricci flow smooths the metric, but it may lead to singularities. It was used in a genius
way by Perelman \cite{Perelman} for solving the Poincar\'e conjecture. Moreover, it was used by Brendle and Schoen \cite{Brendle-Schoen} for
solving the differentiable sphere theorem.

Also, Ricci flow can be applied to discrete geometry, which includes complex networks and weighted finite graphs. To be more specific, we assume that
$G=(V,E,\bf{w})$ is a connected weighted finite graph. Here $V=\{z_1,z_2,\cdots,z_n\}$ denotes the vertex set, $E=\{e_1,e_2,\cdots,e_m\}$ represents the edge set and $\mathbf{w}=(w_{e_1},w_{e_2},\cdots,w_{e_m})$ stands for the weights on edges. In 2009, Ollivier \cite{Ollivier-1} suggested a Ricci flow
\begin{equation}\label{Olloivier-flow}
\texorpdfstring{w_e^\prime(t)=-\kappa_e^{\alpha}(t)w_e(t)}{},
\end{equation}
where, for each $e\in E$, $\kappa_e^\alpha$ is its Ollivier's Ricci curvature, and $\alpha\in[0,1]$ is a parameter. Ni et al \cite{Ni-Lin} found that discrete Ricci flow can be applied to detect community structures in networks, similar to classical Ricci flow on manifolds.

Community detection is an important research area in network analysis, aiming to identify clusters or groups of nodes that are more densely connected internally than with the rest of the network. This concept is crucial in various fields, including sociology \cite{Scott, Wasserman }, biology \cite{Girvan M,Bhowmick}, and computer science\cite{Tauro S L}. Lots of algorithms \cite{Yang-Algesheimer,Fortunato,Newman M E J,Leskovec,Clauset-Newman-Moore,Peel,Girvan M} have been developed to detect and separate communities. Most of these algorithms focus on identifying dense clusters in a graph, using randomized approaches like label propagation or random walks, optimizing centrality measures such as betweenness centrality, or considering measures like modularity.
Ricci curvature, introduced by Forman \cite{Forman}, Ollivier \cite{Ollivier-1} and Lin-Lu-Yau \cite{Lin-Lu-Yau}, acknowledged as an essential tool in graph theory, facilitates deeper insights into network structures and provides an innovative approach to community detection.
In 2019, based on (\ref{Olloivier-flow}) and a surgery process, Ni et al \cite{Ni-Lin} established a very effective method of community detection. In 2021, Bai et al \cite{Bai-2} studied a star coupling Ricci curvature based on Ollivier's Ricci curvature. In 2022, this community detection method was extended by Lai et al \cite{Lai X} to a normalized Ricci flow based on \cite{Lin-Lu-Yau, Bai-2}. Experiments in \cite{Lai X} had shown that their methods were still effective. Very recently, the authors in the current paper \cite{M-Y1} proposed a modified Ricci flow and a quasi-normalized Ricci flow for arbitrary weighted graphs. Notably, the results on global existence and uniqueness of solutions in \cite{M-Y1} do not rely on the exit condition suggested by Bai et al. \cite{Bai-Lin}.

Although literatures \cite{Ni-Lin} and \cite{Lai X} have shown great success in their applications, they lack theoretical support. In order to make up for this shortcoming, recently, concerning Lin-Lu-Yau's Ricci curvature $\kappa_e$, Bai et al \cite{Bai-Lin} obtained the existence and uniqueness of global solutions to the Ricci flow
\begin{equation}\label{flow}
w_e^\prime(t)=-\kappa_e(t)w_e(t)
\end{equation}
and the normalized Ricci flow
\begin{equation}\label{Lai-1-1}
w_e^\prime(t)=-\kappa_e(t)w_e(t) +w_e(t)\sum_{\tau\in E}\kappa_\tau(t) w_\tau(t),
\end{equation}
under an exit condition: if $w_e>\rho_e=\inf_{\gamma}\sum_{\tau\in\gamma}w_\tau$, then the edge $e$ is deleted, where
$e=xy\in E$ and $\gamma$ is a path connecting $x$ and $y$.  Notice that in all of the above mentioned articles,
the authors ask for surgery while extending solutions.

In \cite{M-Y1}, we find that if $w_e$ is replaced by $\rho_e$ on the right hand sides of (\ref{flow}) and
analog of (\ref{Lai-1-1}), namely
 \begin{equation}\label{flow-rho}
w_e^\prime(t)=-\kappa_e(t)\rho_e(t)
\end{equation}
and
\begin{equation}\label{flow-rhoN}
w_e^\prime(t)=-\kappa_e(t)\rho_e(t) +\frac{\sum_{\tau\in E}\kappa_\tau(t)w_\tau(t)}{\sum_\tau w_\tau} \rho_e(t),
\end{equation}
then
both (\ref{flow-rho}) and (\ref{flow-rhoN}) have unique global solutions $w_e(t)$ for each $e\in E$ and all $t\in[0,+\infty)$.
In this setting, we removed the exit condition in \cite{Bai-Lin}, and did not require surgery in the process of flowing. Meanwhile,
experiments show that (\ref{flow-rho}) and (\ref{flow-rhoN}) have the same impressive performance as
(\ref{flow}) and (\ref{Lai-1-1}), and have better modularity.

Currently, there is no existence result for solution to Ollivier's Ricci flow (\ref{Olloivier-flow}). To ensure the global existence of the solution, we can modify the equation (\ref{Olloivier-flow}) by replacing the term $\kappa_e^\alpha w_e$ with $\kappa_e^\alpha \rho_e$. Then we prove the existence and uniqueness of global solutions to the modified Ricci flow. In fact, we shall consider a more general evolution of weights on a connected weighted finite graph $(V,E,\bf{w})$, which reads as
\begin{equation}\label{general-flow}w_e^\prime(t)=W(\mu(x,\cdot),\mu(y,\cdot))-\rho_e,\end{equation}
where $e=xy\in E$ and $W(\mu(x,\cdot),\mu(y,\cdot))$ denotes the Wasserstein distance between two probability measures $\mu(x,\cdot)$ and $\mu(y,\cdot)$. One can easily see that if $\mu(x,z)=\mu_x^\alpha(z)$ is an $\alpha$-lazy one-step random walk, then (\ref{general-flow}) reduces
to $w_e^\prime=-\kappa_e^\alpha\rho_e$. Through a large number of experiments, we conclude that (\ref{general-flow}) is effective in community detection. We look at how key factors, like $\alpha$ and the number of iterations, affect the results. After being tested on real data and compared with existing methods, our approach demonstrates strong performance, especially in measuring modularity.

Before ending this introduction, we mention another interesting curvature flow, the Bakry-\'Emery curvature flow, studied by Cushing et al. Interested readers are referred to  \cite{Cushing-1,Cushing-2} for details.
The remaining part of this paper is organized as follows:  In Section \ref{Sec 2}, we state our main results for the evolution of weights;
In Section \ref{S3}, we prove the Lipschitz property of the Wasserstein distance; The proof of Theorems \ref{theorem1} and \ref{theorem2}  will be given in Section \ref{Sec 3}; In Section \ref{Examples}, We construct several examples of converging flow; Experiments will be done in Section \ref{experiment}.

\section{Notations and main results}\label{Sec 2}

Let $G=(V,E,\mathbf{w})$ be a connected weighted graph, $V=\{z_1,z_2,\cdots,z_n\}$ denotes the set of all vertices, $E=\{e_1,e_2,\cdots,e_m\}$
denotes the set of all edges, and
$\mathbf{w}=(w_{e_1},w_{e_2},\cdots,w_{e_m})$ denotes the vector of weights on edges.
Clearly, each function $f:G\rightarrow\mathbb{R}$ corresponds to a vector $(f(z_1),f(z_2),\cdots,f(z_n))$ in $\mathbb{R}^n$.
Nevertheless, when
the weights $\mathbf{w}$ change, $f$ can be regarded as a vector-valued function, namely
\begin{align*}
{\bf f}:\,\mathbb{R}^m_+\rightarrow{}& \mathbb{R}^n\\[1.2ex]
{\bf{w}}\mapsto{}& (f(z_1),f(z_2),\cdots,f(z_n)).
\end{align*}
Here and in the sequel, $\mathbb{R}^m_+=\{\mathbf{y}=(y_1,y_2,\cdots,y_m)\in\mathbb{R}^m:y_i>0,i=1,2,\cdots,m\}$.
A function $f: G\rightarrow \mathbb{R}$ is said to be {\it locally Lipschitz} in $\mathbb{R}^m_+$ with respect to $\mathbf{w}$, if for any fixed domain $\Omega\subset\subset\mathbb{R}^m_+$, there
exists a constant $C$ depending only on $\Omega$ such that for any vertex $u$,
\begin{equation}\label{Lip}
|f(u)-\widetilde{f}(u)|\leq C|w-\widetilde{w}|,\quad\forall \mathbf{w}, \widetilde{\mathbf{w}}\in\Omega,
\end{equation}
where $\widetilde{f}$ represents the function obtained by replacing $\mathbf{w}$ with
$\widetilde{\mathbf{w}}$ in the expression of $f$.

Throughout this paper, we use the distance between two vertices $x$ and $z$ defined by
\begin{equation}\label{distance-uasual}d(x,z)=\inf_{\gamma}\sum_{\tau\in \gamma}w_\tau,\end{equation}
where the infimum is taken over all paths $\gamma$ connecting $x$ and $z$. Obviously, the function $f=d(x,\cdot)$ is considered
not only  as a vector $(d(x,z_1),d(x,z_2),\cdots,d(x,z_n))$, but also a vector-valued function
${\bf{w}}\mapsto (d(x,z_1),d(x,z_2),\cdots,d(x,z_n))$, in particular, each
$d(x,z_j)$ is a function of $\mathbf{w}$. Given $\mathbf{w}$ and $\widetilde{\mathbf{w}}$ in $\mathbb{R}^m_+$. Let
$d$ and $\widetilde{d}$ be two distance functions determined by $\mathbf{w}$ and $\widetilde{\mathbf{w}}$ respectively.
Then, by (\cite{M-Y1}, Lemma 2.1), for any two fixed vertices $x$ and $z$,
\begin{equation}\label{d-Lip}
|d(x,z)-\widetilde{d}(x,z)|\leq \sqrt{m}|\mathbf{w}-\widetilde{\mathbf{w}}|.
\end{equation}
Thus the distance function $d$ is Lipschitz in $\mathbf{w}$.

 We are now defining a set of functions
\begin{equation}\label{probability}
\mathscr{M}=\left\{\mu:V\times V\rightarrow [0,1]: \sum_{z\in V}\mu(x,z)=1,\,\forall x\in V\right\}.
\end{equation}
Note that for any fixed $x\in V$, each $\mu(x,\cdot)$ is a probability measure. In practical problems, probability measures involving the weights $\mathbf{w}$ are more meaningful.
For example, given any number $\alpha\in[0,1]$, the $\alpha$-lazy one-step random walk reads as
\begin{equation}\label{alpha-lazy}\mu_x^\alpha(z)=\left\{\begin{array}{lll}
\alpha&{\rm if}& z=x\\[1.2ex]
(1-\alpha)\frac{w_{xz}}{\sum_{y\sim x}w_{xy}}&{\rm if}& z\sim x\\[1.2ex]
0&{\rm if}& {\rm otherwise}.
\end{array}\right.\end{equation}
Obviously $\mu(x,\cdot)=\mu_x^\alpha(\cdot)$ belongs to $\mathscr{M}$. Denote the one-step neighbor of $x$ by $N_x=\{y: y\sim x\}$, and
the two-step neighbor of $x$ by $N_{2,x}=\{z\not=x: z\sim y\,\,
{\rm for\,\,some}\,\, y\in N_x\}$.
Similarly, an $\alpha$-lazy two-step random walk is represented by
\begin{equation}\label{2-lazy}\mu_{2,x}^\alpha(z)=\left\{\begin{array}{lll}
\alpha&{\rm if}& z=x\\[1.2ex]
(1-\alpha)\alpha\frac{w_{xz}}{\sum_{y\in N_x}w_{xy}}&{\rm if}& z\in N_x\\[1.2ex]
(1-\alpha)^2\sum_{y\in N_x}\frac{w_{xy}}{\sum_{u\in N_x}w_{xu}}\frac{w_{yz}}{\sum_{v\in N_{2,x}\setminus  N_x}
w_{yv}}&{\rm if}&z\in N_{2,x}\setminus N_x \\[1.2ex]
0&{\rm if}& z\not\in \{x\}\cup N_x\cup N_{2,x}.
\end{array}\right.\end{equation}
One can easily check that $\mu(x,\cdot)=\mu_{2,x}^\alpha(\cdot)\in \mathscr{M}$ for all $x\in V$ and all $\alpha\in[0,1]$. Clearly, $\mu_x^\alpha$ and $\mu_{2,x}^\alpha$ are also
two functions of $\mathbf{w}$, and both of them are locally Lipschitz in $\mathbb{R}^m_+$ with respect to $\mathbf{w}$.

Let $\mu_1$ and $\mu_2$ be two probability measures. A coupling between $\mu_1$ and $\mu_2$ is defined as
a map $A:V\times V\rightarrow[0,1]$ satisfying for all $x, y\in V$,
$$\sum_{u\in V}A(x,u)=\mu_1(x),\quad\sum_{u\in V}A(u,y)=\mu_2(y).$$
  While the Wasserstein distance
between  $\mu_1$ and $\mu_2$ reads
$$W(\mu_1,\mu_2)=\inf_A\sum_{u,v\in V}A(u,v)d(u,v),$$
where the infimum is taken over all couplings between  $\mu_1$ and $\mu_2$.

Assuming that $\mu\in\mathscr{M}$, we consider an evolution of weights $\mathbf{w}=(w_{e_1},w_{e_2},\cdots,w_{e_m})$ according to
\begin{equation}\label{evolution}
\left\{\begin{array}{lll}
w_{e_i}^\prime(t)=W(\mu(x_i,\cdot),\mu(y_i,\cdot))-d(x_i,y_i)\\[1.2ex]
e_i=x_iy_i\in E\\[1.2ex]
w_{e_i}(0)=w_{0,i},\,\,
i=1,2,\cdots,m.
\end{array}\right.
\end{equation}

Our first result reads as follows.

\begin{theorem}\label{theorem1}
Let $G=(V,E,\mathbf{w}_0)$ be a connected weighted finite graph, where
$V$ is the vertex set, $E=\{e_1,e_2,\cdots,e_m\}$ is the edge set, and $\mathbf{w}_0=
(w_{0,1},w_{0,2},\cdots,w_{0,m})\in\mathbb{R}^m_+$ is an arbitrary weight on $E$. If $\mu\in\mathscr{M}$ satisfies that for each vertex
$x\in V$,
$\mu(x,\cdot)$ is locally Lipschitz with respect to $\mathbf{w}_0\in
\mathbb{R}^m_+$,
then the flow (\ref{evolution}) has a unique
solution $\mathbf{w}(t)=(w_{e_1}(t),w_{e_2}(t),\cdots,w_{e_m}(t))$ for $t\in[0,+\infty)$.
\end{theorem}

Analogous to \cite{M-Y1}, we consider a quasi-normalized flow
\begin{equation}\label{evolution-norm}
\left\{\begin{array}{lll}
w_{e_i}^\prime(t)=W\left(\mu(x_i,\cdot),\mu(y_i,\cdot)\right)-d(x_i,y_i)-\frac{\sum_{j=1}^m\left(W(\mu(x_j,\cdot),\mu(y_j,\cdot))-d(x_j,y_j)\right)}
{\sum_{j=1}^mw_{e_j}} d(x_i,y_i)\\[1.2ex]
e_i=x_iy_i\in E\\[1.2ex]
w_{e_i}(0)=w_{0,i},\,\,
i=1,2,\cdots,m.
\end{array}\right.
\end{equation}

\begin{theorem}\label{theorem2}
Under the same assumptions as in Theorem \ref{theorem1}, the flow (\ref{evolution-norm}) has a unique
solution  $\mathbf{w}(t)=(w_{e_1}(t),w_{e_2}(t),\cdots,w_{e_m}(t))$ for $t\in[0,+\infty)$.
\end{theorem}

As we mentioned above, $\alpha$-lazy random walks $\mu_x^\alpha(y)$ and $\mu_{2,x}^\alpha(y)$ are locally Lipschitz in $\mathbf{w}_0\in\mathbb{R}^m_+$. As a consequence, both $\mu(x,y)=\mu_x^\alpha(y)$ and $\nu(x,y)=\mu_{2,x}^\alpha(y)$ are appropriate candidates for $\mu(\cdot,\cdot)$ satisfying the assumptions in Theorems \ref{theorem1} and \ref{theorem2}.
Both the proofs of Theorems \ref{theorem1} and \ref{theorem2} are based on the following lines: Consider the differential system
$$\left\{\begin{array}{lll}
\mathbf{w}^\prime(t)=\mathbf{f}({\mathbf{w}}(t))\\[1.2ex]
\mathbf{w}(0)=\mathbf{w}_0,
\end{array}\right.$$
where ${\bf{f}}(\mathbf{w}(t))=({f}_1(\mathbf{w}(t)),{f}_2(\mathbf{w}(t)),\cdots,{f}_m(\mathbf{w}(t)))$ represents the right hand side of the system (\ref{evolution}) or (\ref{evolution-norm}).
We first prove the local Lipschitz continuity of $\mathbf{w}$. Together with the ODE theory (\cite{Wang-Zhou-Zhu-Wang}, Chapter 6), this implies the local existence and uniqueness of solutions.
Next, we prove there exists a constant $C$ such that
$$-Cw_i\leq {f}_i(\mathbf{w})\leq C\sum_{j=1}^mw_j,\quad\forall i=1,2,\cdots,m.$$
This leads to long time existence of solutions. More details are left to Section \ref{Sec 3}.\\

It is proved by Li-M\"unch \cite{Li-Munch} that the discrete-time Ollivier Ricci curvature flow
$d_{n+1}=(1-\alpha\kappa_{d_n})d_n$ converges to a constant curvature metric, which confirms a conjecture in \cite{Ni-Lin}. However,
addressing the convergence of solutions to systems (\ref{evolution}) or (\ref{evolution-norm}) poses a significant challenge, and we do
not know how to obtain general convergence results except for a few special cases involving particular graphs (see Section \ref{Examples} below). Even so, the global existence of solution is sufficient for practical applications, such as in community detection problems, which will be discussed in  Section \ref{experiment}.

\section{A key lemma}\label{S3}

In this section, we prove a key lemma to be used later. Using the same notations as in Section \ref{Sec 2}. we have the following:

\begin{lemma}\label{Tansport-Lip}
If $\mu\in\mathscr{M}$ and, for each vertex $x$, the probability measure $\mu(x,\cdot)$ is locally Lipschitz in $\mathbf{w}=(w_{e_1},w_{e_2},\cdots,w_{e_m})\in\mathbb{R}^m_+$, then
the Wasserstein distance $W\left(\mu(x,\cdot),\mu(y,\cdot)\right)$ is also locally Lipschitz in $\mathbf{w}\in\mathbb{R}^m_+$.
\end{lemma}

\begin{proof} Given two vertices $x$, $y$ and two weights ${\bf{w}}=(w_{e_1},\cdots,w_{e_m})$, $\widetilde{{\bf{w}}}=(\widetilde{w}_{e_1},\cdots,\widetilde{w}_{e_m})$ in $\mathbb{R}^m_+$. Assume that Wasserstein distances  $W(\mu(x,\cdot),\mu(y,\cdot))$ and $W(\widetilde{\mu}(x,\cdot),\widetilde{\mu}(y,\cdot))$
are determined by ${\bf{w}}$ and $\widetilde{{\bf{w}}}$ respectively, as well as $\mu$ and $\widetilde{\mu}$.  With no loss of generality, we  assume there exist some constants $\Lambda>0$ and $\delta>0$ such that for each $i=1,2,\cdots,m$,
\begin{equation}\label{hypo}\Lambda^{-1}\leq w_{e_i}\leq\Lambda,\,\,
\Lambda^{-1}\leq \widetilde{w}_{e_i}\leq\Lambda,\,\,
|w_{e_i}-\widetilde{w}_{e_i}|\leq \delta. \end{equation}
By the Kantorovich-Rubinstein duality formula,
\begin{align}
W\left(\mu(x,\cdot),\mu(y,\cdot)\right)&=\sup_{\psi\in{\rm Lip}\,1}\sum_{u\in V} \psi(u)(\mu(x,u)-\mu(y,u)),\label{W1}\\[1.2ex]
W(\widetilde{\mu}(x,\cdot),\widetilde{\mu}(y,\cdot))&=\sup_{\psi\in\widetilde{{\rm Lip}}\,1}\sum_{u\in V} \psi(u)(\widetilde{\mu}(x,u)-\widetilde{\mu}(y, u)),\label{W2}
\end{align}
where
\begin{gather*}
{\rm Lip}\,1=\left\{f:V\rightarrow\mathbb{R}: |f(u)-f(v)|\leq d(u,v),\,\forall u,v\in V\right\},\\[1.2ex]
\widetilde{{\rm Lip}}\,1=\{f:V\rightarrow\mathbb{R}: |f(u)-f(v)|\leq \widetilde{d}(u,v),\,\forall u,v\in V\},
\end{gather*}
the distance functions $d$ and $\widetilde{d}$ are determined by $\mathbf{w}$ and $\widetilde{\mathbf{w}}$ respectively.
Now we distinguish two cases to proceed. \\

{\bf Case 1}. $W(\mu(x,\cdot),\mu(y,\cdot))\geq W(\widetilde{\mu}(x,\cdot),\widetilde{\mu}(y,\cdot))$.\\

In view of (\ref{W1}), there exists some $f\in{\rm Lip}\,1$ such that
$$\sum_{u\in V} f(u)(\mu(x,u)-\mu(y,u))=\sup_{\psi\in{\rm Lip}\,1}\sum_{u\in V} \psi(u)(\mu(x,u)-\mu(y,u)).$$
By (\cite{M-Y1}, Lemma 2.1), we have for all vertices $u$ and $v$,
$$|d(u,v)-\widetilde{d}(u,v)|\leq\sqrt{m}|\mathbf{w}-\widetilde{\mathbf{w}}|.$$
It then follows that
\[
\frac{d(u,v)}{\widetilde{d}(u,v)}=1+\frac{d(u,v)-\widetilde{d}(u,v)}{\widetilde{d}(u,v)}
\leq1+\sqrt{m}\Lambda|\mathbf{w}-\widetilde{\mathbf{w}}|.
\]
Set
$$\widetilde{f}(u)=\frac{f(u)}{1+\sqrt{m}\Lambda|\mathbf{w}-\widetilde{\mathbf{w}}|},\quad\forall u\in V.$$
Since
$$|\widetilde{f}(u)-\widetilde{f}(v)|=\frac{|f(u)-f(v)|}{1+\sqrt{m}\Lambda|\mathbf{w}-\widetilde{\mathbf{w}}|}\leq \widetilde{d}(u,v)$$
for all $u,v\in V$, we have $\widetilde{f}\in \widetilde{{\rm Lip}}\,1$.
By our assumption $\mu(x,\cdot)$ is locally Lipschitz in $\mathbb{R}^m_+$ with respect to $\mathbf{w}$,
there exists a constant $C_0>0$, depending only on $\Lambda$ and $\delta$, such that for all $\mathbf{w}$ and $\widetilde{\mathbf{w}}$ satisfying (\ref{hypo}),
$$|\mu(x,u)-\widetilde{\mu}(x,u)|\leq C_0|\mathbf{w}-\widetilde{\mathbf{w}}|,\quad \forall x,u\in V.$$
As a consequence, we obtain
\begin{align*}
&W(\mu(x,\cdot),\mu(y,\cdot))-W(\widetilde{\mu}(x,\cdot),\widetilde{\mu}(y,\cdot))\\[1.2ex]
={}&\sum_{u\in V}f(u)(\mu(x,u)-\mu(y,u))
-\sup_{\psi\in \widetilde{{\rm Lip}}\,1}\sum_{u\in V}\psi(u)(\widetilde{\mu}(x,u)-\widetilde{\mu}(y,u))\\[1.2ex]
\leq{}&\sum_{u\in V}f(u)(\mu(x,u)-\mu(y,u))-\sum_{u\in V}\widetilde{f}(u)(\widetilde{\mu}(x,u)-\widetilde{\mu}(y,u))\\[1.2ex]
\leq{}& \sum_{u\in V}|f(u)|(|\mu(x,u)-\widetilde{\mu}(x,u)|+|\mu(y,u)-\widetilde{\mu}(y,u)|)\\[1.2ex]
&\quad+\sum_{u\in V} |f(u)-\widetilde{f}(u)|(\widetilde{\mu}(x,u)+\widetilde{\mu}(y,u))\\[1.2ex]
\leq{}& 2n(C_0+\sqrt{m}\Lambda)\|f\|_{L^\infty(V)}|\mathbf{w}-\widetilde{\mathbf{w}}|,
\end{align*}
where $n$ is the number of all vertices of $V$.\\

{\bf Case 2}. $W(\mu(x,\cdot),\mu(y,\cdot))< W(\widetilde{\mu}(x,\cdot),\widetilde{\mu}(y,\cdot))$.\\

By swapping the positions of $\mathbf{w}$ and $\widetilde{\mathbf{w}}$, we have by the same argument as in Case 1,
$$W(\widetilde{\mu}(x,\cdot),\widetilde{\mu}(y,\cdot))-W(\mu(x,\cdot),\mu(y,\cdot))\leq 2n(C_0+\sqrt{m}\Lambda)\|f\|_{L^\infty(V)}|\mathbf{w}-\widetilde{\mathbf{w}}|.$$

Combining Cases 1 and  2, we complete the proof of the lemma.
\end{proof}

\section{Proof of Theorems \ref{theorem1} and \ref{theorem2}}\label{Sec 3}

In this section, we prove Theorems \ref{theorem1} and \ref{theorem2} by using the ODE theory. \\

{\it Proof of Theorem \ref{theorem1}.}
 We divide the proof into two parts.\\

{\bf Part 1.} {\it Short time existence.}\\

Set $w_i=w_{e_i}$, $i=1,2,\cdots,m$. Given any vector
${\bf{w}}_0=(w_{0,1},w_{0,2},\cdots,w_{0,m})\in\mathbb{R}^m_+$. In fact, the evolution of weights (\ref{evolution})
is the ordinary differential system
 \begin{equation}\label{equiv}
 \left\{\begin{array}{lll}
 {\bf{w}}^\prime(t)={\bf f}({\bf{w}}(t))\\[1.5ex]
 {{\bf{w}}(0)}={{\bf{w}}_0},
 \end{array}\right.
 \end{equation}
where ${\bf w}=(w_1,w_2,\cdots,w_m)\in\mathbb{R}^m_+$ and
 ${\mathbf{f}}=(f_1,f_2,\cdots,f_m)$ is a map represented by
\begin{align*}
{\bf f}:\mathbb{R}^m_+\rightarrow{}& \mathbb{R}^m\\
{\bf{w}}\mapsto{}& (f_1(\mathbf{w}),\cdots,f_m(\mathbf{w})),
\end{align*}
where $f_i(\mathbf{w})=W(\mu(x_i,\cdot),\mu(y_i,\cdot))-d(x_i,y_i)$ and $e_i=x_iy_i$, $i=1,2,\cdots,m$. From Lemma \ref{Tansport-Lip} and
(\cite{M-Y1}, Lemma 2.1), we know that all $f_i(\mathbf{w})$, $i=1,2,\cdots, m$, are locally Lipschitz in
$\mathbb{R}^m_+$, and so is ${\bf f}({\bf w})$. By the ODE theory (\cite{Wang-Zhou-Zhu-Wang}, Chapter 6), there exists a constant $T>0$ such that the ordinary differential system (\ref{equiv}) has a unique
solution ${\bf{w}}(t)$ on $[0,T]$. \\

{\bf Part 2.} {\it Long time existence.}\\

According to the conclusion of the first part, we may define
$$T^\ast=\sup\{T>0: (\ref{equiv})\, {\rm has\,a\,unique\,solution\,on\,}[0,T]\}.$$
If $T^\ast<+\infty$, then (\ref{equiv}) has a unique solution ${\bf w}(t)$ on the time interval
$[0,T^\ast)$; moreover, according to the ODE theory (\cite{Wang-Zhou-Zhu-Wang}, Chapter 6), we have either
\begin{equation}\label{tend-1}\liminf_{t\rightarrow T^\ast}\phi(t)= 0\end{equation}
or
\begin{equation}\label{tend-2}\limsup_{t\rightarrow T^\ast}\Phi(t)= +\infty,\end{equation}
where
$$\phi(t)=\min\{w_{1}(t),w_{2}(t),\cdots,w_{m}(t)\}\quad\text{and}\quad\Phi(t)=\max\{w_{1}(t),w_{2}(t),\cdots,w_{m}(t)\}.$$

Let $e=xy\in E$ be fixed. Since
\begin{equation}\label{geq-1}W(\mu(x,\cdot),\mu(y,\cdot))-\rho_e(t)\geq -\rho_e(t)\geq -w_e(t),\end{equation}
this together with (\ref{equiv}) implies
$$w_e^\prime(t)=W(\mu(x,\cdot),\mu(y,\cdot))-\rho_e(t)\geq -w_e(t).$$
Thus we have
\begin{equation}\label{lower-bd}
w_e(t)\geq w_e(0)e^{-T^\ast},\quad\forall t\in[0,T^\ast).
\end{equation}
Noting also that each coupling $A$ between $\mu(x,\cdot)$ and $\mu(y,\cdot)$ satisfies $A(u,v)\in [0,1]$ and
$\sum_{u,v\in V}A(u,v)=1$, and that $d(u,v)\leq \sum_{\tau\in E}w_{\tau}$ for all vertices $u,v$, we obtain
\begin{align}
W(\mu(x,\cdot),\mu(y,\cdot))={}&\inf_{B}\sum_{u,v\in V}B(u,v)d(u,v)\nonumber\\
\leq{}& \sum_{u,v\in V}A(u,v)d(u,v)\nonumber\\
\leq{}&\left(\sum_{u,v\in V}A(u,v)\right)\left(\sum_{\tau\in E}w_{\tau}\right)\nonumber\\
={}&\sum_{\tau\in E}w_{\tau},\label{leq-1}
\end{align}
where, in the first equality, $B$ is taken from all couplings between probability measures $\mu(x,\cdot)$ and $\mu(y,\cdot)$.
In view of (\ref{equiv}), it follows that
\begin{align*}
\frac{d}{dt}\sum_{\tau\in E}w_\tau(t)={}&\sum_{\tau=xy\,\in\, E}\left(W(\mu(x,\cdot),\mu(y,\cdot))-d(x,y)\right)\\
\leq{}& \sum_{\tau=xy\,\in\, E}W(\mu(x,\cdot),\mu(y,\cdot))\\
\leq{}&m\sum_{\tau\in E}w_\tau(t).
\end{align*}
Then integration by parts gives
\begin{equation}\label{sum-up}
\sum_{\tau\in E}w_\tau(t)\leq e^{mT^\ast}\sum_{\tau\in E}w_\tau(0),\quad\forall t\in[0,T^\ast).
\end{equation}
Combining (\ref{lower-bd}) and (\ref{sum-up}), we have
$$\phi(0)e^{-T^\ast}\leq \phi(t)\leq \Phi(t)\leq \sum_{\tau\in E}w_\tau(t)\leq e^{mT^\ast}\sum_{\tau\in E}w_\tau(0)$$
for all $t\in [0,T^\ast)$, which contradicts (\ref{tend-1}) and (\ref{tend-2}).
 Therefore $T^\ast=+\infty$, and thus (\ref{equiv}) has a unique solution $\mathbf{w}(t)$ on $[0,+\infty)$. $\hfill\Box$\\

{\it Proof of Theorem \ref{theorem2}.}

Given ${\bf{w}}=(w_{e_1},w_{e_2},\cdots,w_{e_m})\in \mathbb{R}^m_+$. Let
${\bf{g}}:\mathbb{R}^m_+\rightarrow \mathbb{R}^m$ be a vector-valued function written as ${\bf{g}}({\bf{w}})=(g_1({\bf{w}}),g_2({\bf{w}}),
\cdots,g_m({\bf{w}}))$, where
$${{g}_i}({\bf{w}})=W\left(\mu(x_i,\cdot),\mu(y_i,\cdot)\right)-\rho_{e_i}-\frac{\sum_{j=1}^m\left(W(\mu(x_j,\cdot),\mu(y_j,\cdot))-
\rho_{e_j}\right)}{\sum_{j=1}^mw_{e_j}} \rho_{e_i}$$
and $e_i=x_iy_i$.
From Lemma \ref{Tansport-Lip} and (\cite{M-Y1}, Lemma 2.1), we conclude that ${\bf{g}}$ is locally Lipschitz in $\mathbb{R}^m_+$, namely
if $\Omega\subset\mathbb{R}^m_+$ satisfies $\overline{\Omega}\subset\mathbb{R}^m_+$, then there exists a constant $C$ depending only on
$\Omega$ such that
$$|{\bf{g}}({\bf{w}})-{\bf{g}}({\widetilde{\bf{w}}})|\leq C|{\bf{w}}-{\widetilde{\bf{w}}}|,\quad\forall {\bf{w}}, {\widetilde{\bf{w}}}\in\Omega.$$
Hence, according to the ODE theory (\cite{Wang-Zhou-Zhu-Wang}, Chapter 6), there exists some $T>0$ such that the flow (\ref{evolution-norm}) has a unique solution ${\bf{w}}(t)$ on $[0,T]$. Let
$$T^\ast=\sup\{T>0: (\ref{evolution-norm})\, {\rm has\,a\,unique\,solution\,on\,}[0,T]\}.$$
If $T^\ast<+\infty$, then (\ref{evolution-norm}) has a unique solution ${\bf w}(t)=(w_{e_1}(t),w_{e_2}(t),\cdots,w_{e_m}(t))$ for all
$t\in[0,T^\ast)$. However, the ODE theory (\cite{Wang-Zhou-Zhu-Wang}, Chapter 6) implies either
\begin{equation}\label{tend-3}\liminf_{t\rightarrow T^\ast}\min\{w_{e_1}(t),w_{e_2}(t),\cdots,w_{e_m}(t)\}= 0\end{equation}
or
\begin{equation}\label{tend-4}\limsup_{t\rightarrow T^\ast}\max\{w_{e_1}(t),w_{e_2}(t),\cdots,w_{e_m}(t)\}= +\infty.\end{equation}
In view of (\ref{geq-1}) and (\ref{leq-1}), we obtain for each $i$,
$$-(m+1)w_{e_i}\leq g_i(\mathbf{w})\leq w_{e_i}+\sum_{j=1}^mw_{e_j}.$$
This together with (\ref{evolution-norm}) gives for each $i$,
\begin{equation}\label{ineq-1}-(m+1)w_{e_i}(t)\leq w_{e_i}^\prime(t)\leq w_{e_i}(t)+\sum_{j=1}^mw_{e_j}(t)\end{equation}
and
\begin{equation}\label{ineq-2}\frac{d}{dt}\sum_{i=1}^mw_{e_i}\leq (m+1)\sum_{i=1}^mw_{e_i}.\end{equation}
From the left side of (\ref{ineq-1}),  we conclude that $w_{e_i}(t)\geq w_{e_i}(0)e^{-(m+1)T^\ast}$ for all  $t\in[0,T^\ast)$ and each $i$, contradicting (\ref{tend-3}). While from (\ref{ineq-2}), we conclude that
$$\sum_{\tau\in E}w_\tau(t)\leq e^{(m+1)T^\ast}\sum_{\tau\in E}w_\tau(0)$$
for all $t\in[0,T^\ast)$, contradicting (\ref{tend-4}).
 Therefore $T^\ast=+\infty$, as we desired. $\hfill\Box$

\section{Examples}\label{Examples}

According to Theorems \ref{theorem1} and \ref{theorem2}, we know that both of evolutions of weights
(\ref{evolution}) and (\ref{evolution-norm}) have unique global solutions. However, in general, we do not ensure that these solutions
converge. In this section, we shall only construct examples of convergent solutions to (\ref{evolution}) for special graphs and constant
initial weights, but leave the case (\ref{evolution-norm}) to interested readers.\\

\begin{figure}[H]
    \centering
   \begin{tikzpicture}[x=0.75pt,y=0.75pt,yscale=-1,xscale=1]

\draw    (23.75,265.06) -- (122.55,265.06) ;
\draw  [fill={rgb, 255:red, 0; green, 0; blue, 0 }  ,fill opacity=1 ] (21.28,265.06) .. controls (21.28,263.7) and (22.38,262.59) .. (23.75,262.59) .. controls (25.12,262.59) and (26.23,263.7) .. (26.23,265.06) .. controls (26.23,266.43) and (25.12,267.54) .. (23.75,267.54) .. controls (22.38,267.54) and (21.28,266.43) .. (21.28,265.06) -- cycle ;
\draw  [fill={rgb, 255:red, 0; green, 0; blue, 0 }  ,fill opacity=1 ] (120.08,265.06) .. controls (120.08,263.7) and (121.18,262.59) .. (122.55,262.59) .. controls (123.92,262.59) and (125.03,263.7) .. (125.03,265.06) .. controls (125.03,266.43) and (123.92,267.54) .. (122.55,267.54) .. controls (121.18,267.54) and (120.08,266.43) .. (120.08,265.06) -- cycle ;

\draw (11.56,264) node [anchor=north west][inner sep=0.75pt]   [align=left] {$x$};
\draw (123.75,264.46) node [anchor=north west][inner sep=0.75pt]   [align=left] {$y$};
\draw (62.76,246.2) node [anchor=north west][inner sep=0.75pt]   [align=left] {${w}_0$};

\end{tikzpicture}
    \caption{A line segment}
    \label{example_1}
\end{figure}
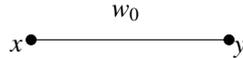

{\it Example 1.} Let $G=(V,E,\mathbf{w}_0)$ be a line segment illustrated in Figure \ref{example_1}, where $V=\{x,y\}$, $E=\{e=xy\}$, and $\mathbf{w}_{0}=w_0$ is the initial weight of the edge $e$.
For any $\alpha\in[0,1]$, the $\alpha$-lazy one-step random walks are written as
$$\mu_x^\alpha(u)=\left\{\begin{array}{lll}
\alpha&{\rm if}& u=x\\
1-\alpha&{\rm if}& u=y,
\end{array}\right.
\quad \mu_y^\alpha(u)=\left\{\begin{array}{lll}
1-\alpha&{\rm if}& u=x\\
\alpha&{\rm if}& u=y.
\end{array}\right.
$$

Write $\mu(x,\cdot)=\mu_x^\alpha(\cdot)$ and $\mu(y,\cdot)=\mu_y^\alpha(\cdot)$.
For $\alpha\in[0,1/2]$, one calculates the Wasserstein distance between $\mu_x$ and $\mu_y$ as
$W(\mu_x^\alpha,\mu_y^\alpha)=(1-2\alpha)w_e$, and the graph distance $d(x,y)=w_{e}$.
As a consequence, the flow (\ref{evolution}) becomes
$$w_e^\prime(t)=W(\mu_x^\alpha,\mu_y^\alpha)-d(x,y)=-2\alpha w_e(t).$$
Thus $w_e(t)=w_{0}e^{-2\alpha t}$ for all $t\in[0,+\infty)$. Similarly, one derives that for $\alpha\in[1/2,1]$, (\ref{evolution}) has the unique
solution $w_e(t)=w_{0}e^{2(\alpha-1)t}$ for all $t\in[0,+\infty)$. In both cases, along the flow (\ref{evolution}), $w_e(t)$ converges to zero exponentially.\\

\begin{figure}[H]
    \centering
   \begin{tikzpicture}[x=0.75pt,y=0.75pt,yscale=-1,xscale=1]

\draw    (23.75,265.06) -- (122.55,265.06) ;
\draw  [fill={rgb, 255:red, 0; green, 0; blue, 0 }  ,fill opacity=1 ] (21.28,265.06) .. controls (21.28,263.7) and (22.38,262.59) .. (23.75,262.59) .. controls (25.12,262.59) and (26.23,263.7) .. (26.23,265.06) .. controls (26.23,266.43) and (25.12,267.54) .. (23.75,267.54) .. controls (22.38,267.54) and (21.28,266.43) .. (21.28,265.06) -- cycle ;
\draw  [fill={rgb, 255:red, 0; green, 0; blue, 0 }  ,fill opacity=1 ] (120.08,265.06) .. controls (120.08,263.7) and (121.18,262.59) .. (122.55,262.59) .. controls (123.92,262.59) and (125.03,263.7) .. (125.03,265.06) .. controls (125.03,266.43) and (123.92,267.54) .. (122.55,267.54) .. controls (121.18,267.54) and (120.08,266.43) .. (120.08,265.06) -- cycle ;
\draw    (122.55,265.06) -- (221.35,265.06) ;
\draw  [fill={rgb, 255:red, 0; green, 0; blue, 0 }  ,fill opacity=1 ] (218.88,265.06) .. controls (218.88,263.7) and (219.98,262.59) .. (221.35,262.59) .. controls (222.72,262.59) and (223.83,263.7) .. (223.83,265.06) .. controls (223.83,266.43) and (222.72,267.54) .. (221.35,267.54) .. controls (219.98,267.54) and (218.88,266.43) .. (218.88,265.06) -- cycle ;

\draw (10.56,266) node [anchor=north west][inner sep=0.75pt]   [align=left] {$x$};
\draw (118.75,268.46) node [anchor=north west][inner sep=0.75pt]   [align=left] {$y$};
\draw (62.76,246.2) node [anchor=north west][inner sep=0.75pt]   [align=left] {${w}_0$};
\draw (223.35,265.59) node [anchor=north west][inner sep=0.75pt]   [align=left] {$z$};
\draw (166.26,247.5) node [anchor=north west][inner sep=0.75pt]   [align=left] {${w}_0$};

\end{tikzpicture}

    \caption{A path graph G with length 2}
    \label{example_2}
\end{figure}
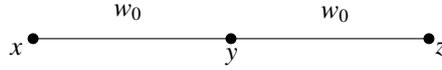

{\it Example 2.} Let $G=(V,E,\mathbf{w}_0)$ be a three-point path in Figure \ref{example_2}, where $V=\{x,y,z\}$, $E=\{xy,yz\}$, and $\mathbf{w}_0=(w_{0},w_{0})$ is the initial weight.
Let $\alpha\in[0,1]$, $\mu(x,\cdot)=\mu_x^\alpha(\cdot)$ and $\mu(y,\cdot)=\mu_y^\alpha(\cdot)$, where
$$\mu_x^\alpha(u)=\left\{\begin{array}{lll}
\alpha&{\rm if}& u=x\\
1-\alpha&{\rm if}& u=y,
\end{array}\right.
\quad   \mu_z^\alpha(u)=\left\{\begin{array}{lll}
\alpha&{\rm if}& u=z\\
1-\alpha&{\rm if}& u=y,
\end{array}\right.
$$ and
$$\mu_y^\alpha(u)=\left\{\begin{array}{lll}
(1-\alpha)\frac{w_{xy}}{w_{xy}+w_{yz}}&{\rm if}& u=x\\
\alpha&{\rm if}& u=y\\
(1-\alpha)\frac{w_{yz}}{w_{xy}+w_{yz}}&{\rm if}&u=z.
\end{array}\right.$$
Assuming that $w_{xy}=w_{yz}$, we have $\mu_y(u)=(1-\alpha)/2$ for $u\in\{x,z\}$, $\mu_y(u)=\alpha$ for $u=y$, and the Wasserstein distance
$W(\mu_x,\mu_y)=W(\mu_y,\mu_z)=(1-2\alpha)w_{xy}$ for $\alpha\in[0,1/3]$, $W(\mu_x,\mu_y)=W(\mu_y,\mu_z)=\alpha w_{xy}$ for $\alpha\in[1/3,1]$,
$d(x,y)=w_{xy}$ and $d(y,z)=w_{yz}$.
Obviously, the initial value problem
$$\left\{\begin{array}{lll}
w_{xy}^\prime(t)=W(\mu_x^\alpha,\mu_y^\alpha)-d(x,y)=-2\alpha w_{xy}\\
w_{yz}^\prime(t)=W(\mu_y,\mu_z)-d(y,z)=-2\alpha w_{yz}\\
w_{xy}(t)=w_{yz}(t)\\
w_{xy}(0)=w_{0}=w_{yz}(0)
\end{array}\right.$$
has a solution
\begin{equation}\label{solu}(w_{xy}(t), w_{yz}(t))=\left\{\begin{array}{lll}(w_0e^{-2\alpha t},w_0e^{-2\alpha t})&{\rm for}&\alpha\in[0,1/3]\\[1.2ex]
(w_0e^{(\alpha-1) t},w_0e^{(\alpha-1) t})&{\rm for}&\alpha\in[1/3,1].\end{array}\right.\end{equation}
From Theorem \ref{theorem1}, we know that (\ref{evolution}) has a unique global solution. Hence (\ref{solu})
is the solution of (\ref{evolution}) with the initial weights $(w_{xy}(0),w_{yz}(0))=(w_0,w_0)$.\\

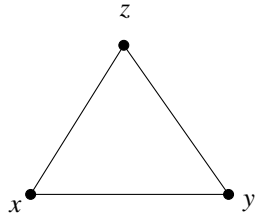
\begin{figure}[H]
    \centering
    \begin{tikzpicture}[x=0.75pt,y=0.75pt,yscale=-1,xscale=1]

\draw    (23.75,265.06) -- (122.55,265.06) ;
\draw  [fill={rgb, 255:red, 0; green, 0; blue, 0 }  ,fill opacity=1 ] (21.28,265.06) .. controls (21.28,263.7) and (22.38,262.59) .. (23.75,262.59) .. controls (25.12,262.59) and (26.23,263.7) .. (26.23,265.06) .. controls (26.23,266.43) and (25.12,267.54) .. (23.75,267.54) .. controls (22.38,267.54) and (21.28,266.43) .. (21.28,265.06) -- cycle ;
\draw  [fill={rgb, 255:red, 0; green, 0; blue, 0 }  ,fill opacity=1 ] (120.08,265.06) .. controls (120.08,263.7) and (121.18,262.59) .. (122.55,262.59) .. controls (123.92,262.59) and (125.03,263.7) .. (125.03,265.06) .. controls (125.03,266.43) and (123.92,267.54) .. (122.55,267.54) .. controls (121.18,267.54) and (120.08,266.43) .. (120.08,265.06) -- cycle ;
\draw    (122.55,265.06) -- (70.27,190.05) ;
\draw  [fill={rgb, 255:red, 0; green, 0; blue, 0 }  ,fill opacity=1 ] (67.79,190.05) .. controls (67.79,188.69) and (68.9,187.58) .. (70.27,187.58) .. controls (71.63,187.58) and (72.74,188.69) .. (72.74,190.05) .. controls (72.74,191.42) and (71.63,192.53) .. (70.27,192.53) .. controls (68.9,192.53) and (67.79,191.42) .. (67.79,190.05) -- cycle ;
\draw    (23.75,265.06) -- (70.27,190.05) ;

\draw (11.23,266.67) node [anchor=north west][inner sep=0.75pt]   [align=left] {$x$};
\draw (128.55,263.59) node [anchor=north west][inner sep=0.75pt]   [align=left] {$y$};
\draw (67.02,168.25) node [anchor=north west][inner sep=0.75pt]   [align=left] {$z$};

\end{tikzpicture}

    \caption{A triangle}
    \label{example_3}
\end{figure}

{\it Example 3.} Let $G=(V,E,\mathbf{w}_0)$ be a triangle in Figure \ref{example_3}, where $V=\{x,y,z\}$, $E=\{xy,yz,zx\}$, $\mathbf{w}_0=(w_0,w_0,w_0)$.
Assume that $\mathbf{w}=(w,w,w)$ is another weight of $G$, under which we have
 $$\mu_x^\alpha(u)=\left\{\begin{array}{lll}
\alpha,& u=x\\
\frac{1-\alpha}{2},& u=y,\,z,
\end{array}\right.\quad \mu_y^\alpha(u)=\left\{\begin{array}{lll}
\frac{1-\alpha}{2},& u=x,\,z\\
\alpha,& u=y.
\end{array}\right. $$
This leads to
$$W(\mu_x^\alpha,\mu_y^\alpha)=\left\{\begin{array}{lll}
\frac{1-3\alpha}{2}w_0,&\alpha\in[0,1/3]\\[1.2ex]
\frac{3\alpha-1}{2}w_0,&\alpha\in[1/3,1].
\end{array}\right.$$
For $\alpha\in[0,1/3]$, the initial value problem
$$\left\{\begin{array}{lll}
w^\prime(t)=-\frac{3\alpha+1}{2}w(t)\\
w(0)=w_0
\end{array}\right.$$
has a solution $w(t)=w_0e^{-\frac{3\alpha+1}{2}t}$. For $\alpha\in[1/3,1]$, the initial value problem
$$\left\{\begin{array}{lll}
w^\prime(t)=\frac{3(\alpha-1)}{2}w(t)\\
w(0)=w_0
\end{array}\right.$$
has a solution $w(t)=w_0e^{\frac{3(\alpha-1)}{2}t}$. Let $\mu(x,\cdot)=\mu_x^\alpha(\cdot)$, $\mu(y,\cdot)=\mu_y^\alpha(\cdot)$
and $\mu(z,\cdot)=\mu_z^\alpha(\cdot)$. Then the unique global solution of the flow (\ref{evolution}) is
$(w(t),w(t),w(t))$.  \\

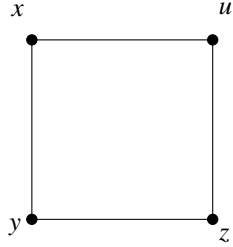
\begin{figure}[H]
    \centering
\begin{tikzpicture}[x=0.75pt,y=0.75pt,yscale=-1,xscale=1]

\draw   (50,171) -- (140.27,171) -- (140.27,261.27) -- (50,261.27) -- cycle ;
\draw  [fill={rgb, 255:red, 0; green, 0; blue, 0 }  ,fill opacity=1 ] (52.57,171) .. controls (52.57,169.58) and (51.42,168.43) .. (50,168.43) .. controls (48.58,168.43) and (47.43,169.58) .. (47.43,171) .. controls (47.43,172.42) and (48.58,173.57) .. (50,173.57) .. controls (51.42,173.57) and (52.57,172.42) .. (52.57,171) -- cycle ;
\draw  [fill={rgb, 255:red, 0; green, 0; blue, 0 }  ,fill opacity=1 ] (142.83,171) .. controls (142.83,169.58) and (141.68,168.43) .. (140.27,168.43) .. controls (138.85,168.43) and (137.7,169.58) .. (137.7,171) .. controls (137.7,172.42) and (138.85,173.57) .. (140.27,173.57) .. controls (141.68,173.57) and (142.83,172.42) .. (142.83,171) -- cycle ;
\draw  [fill={rgb, 255:red, 0; green, 0; blue, 0 }  ,fill opacity=1 ] (52.57,261.27) .. controls (52.57,259.85) and (51.42,258.7) .. (50,258.7) .. controls (48.58,258.7) and (47.43,259.85) .. (47.43,261.27) .. controls (47.43,262.68) and (48.58,263.83) .. (50,263.83) .. controls (51.42,263.83) and (52.57,262.68) .. (52.57,261.27) -- cycle ;
\draw  [fill={rgb, 255:red, 0; green, 0; blue, 0 }  ,fill opacity=1 ] (142.83,261.27) .. controls (142.83,259.85) and (141.68,258.7) .. (140.27,258.7) .. controls (138.85,258.7) and (137.7,259.85) .. (137.7,261.27) .. controls (137.7,262.68) and (138.85,263.83) .. (140.27,263.83) .. controls (141.68,263.83) and (142.83,262.68) .. (142.83,261.27) -- cycle ;

\draw (38,151.67) node [anchor=north west][inner sep=0.75pt]   [align=left] {$x$};
\draw (37.33,259.33) node [anchor=north west][inner sep=0.75pt]   [align=left] {$y$};
\draw (142.27,264.27) node [anchor=north west][inner sep=0.75pt]   [align=left] {$z$};
\draw (142,151) node [anchor=north west][inner sep=0.75pt]   [align=left] {$u$};

\end{tikzpicture}

    \caption{A square}
    \label{example_4}
\end{figure}

{\it Example 4.} $G=(V,E,\mathbf{w}_0)$ is a square in Figure \ref{example_4}, $V=\{x,y,z,u\}$, $E=\{xy,yz,zu,ux\}$, $\mathbf{w}_0=(w_0,w_0,w_0,w_0)$.
For weights $\mathbf{w}=(w,w,w,w)$, we have for $\alpha\in[0,1]$,
$$\mu_x^\alpha(v)=\left\{\begin{array}{lll}
\alpha,& v=x\\
\frac{1-\alpha}{2},& v=y,\,u\\
0,&v=z,
\end{array}\right.\quad \mu_y^\alpha(v)=\left\{\begin{array}{lll}
\alpha,& v=y\\
\frac{1-\alpha}{2},& v=x,\,z\\
0,&v=u
\end{array}\right.$$
and
$$W(\mu_x^\alpha,\mu_y^\alpha)=\left\{\begin{array}{lll}
(1-2\alpha)w,&\alpha\in[0,1/3]\\
\alpha w,&\alpha\in[1/3,1].
\end{array}\right.$$
This together with the differential equation
$$w^\prime(t)=w_{xy}^\prime(t)=W(\mu_x^\alpha,\mu_y^\alpha)-d(x,y)$$
gives
\begin{equation}\label{at}w(t)=\left\{\begin{array}{lll}
w_0e^{-2\alpha t} ,&\alpha\in[0,1/3]\\
w_0e^{(\alpha-1)t} ,&\alpha\in[1/3,1].
\end{array}\right.\end{equation}
Set $\mu(v,\cdot)=\mu_v(\cdot)$ for all $v\in V$. The evolution problem (\ref{evolution}) has a unique  global solution
$\mathbf{w}(t)=(w(t),w(t),w(t),w(t))$ with $w(t)$ given as in (\ref{at}).\\

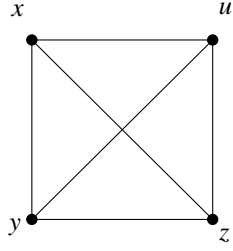
\begin{figure}[H]
    \centering
    \begin{tikzpicture}[x=0.75pt,y=0.75pt,yscale=-1,xscale=1]

\draw   (50,171) -- (140.27,171) -- (140.27,261.27) -- (50,261.27) -- cycle ;
\draw  [fill={rgb, 255:red, 0; green, 0; blue, 0 }  ,fill opacity=1 ] (52.57,171) .. controls (52.57,169.58) and (51.42,168.43) .. (50,168.43) .. controls (48.58,168.43) and (47.43,169.58) .. (47.43,171) .. controls (47.43,172.42) and (48.58,173.57) .. (50,173.57) .. controls (51.42,173.57) and (52.57,172.42) .. (52.57,171) -- cycle ;
\draw  [fill={rgb, 255:red, 0; green, 0; blue, 0 }  ,fill opacity=1 ] (142.83,171) .. controls (142.83,169.58) and (141.68,168.43) .. (140.27,168.43) .. controls (138.85,168.43) and (137.7,169.58) .. (137.7,171) .. controls (137.7,172.42) and (138.85,173.57) .. (140.27,173.57) .. controls (141.68,173.57) and (142.83,172.42) .. (142.83,171) -- cycle ;
\draw  [fill={rgb, 255:red, 0; green, 0; blue, 0 }  ,fill opacity=1 ] (52.57,261.27) .. controls (52.57,259.85) and (51.42,258.7) .. (50,258.7) .. controls (48.58,258.7) and (47.43,259.85) .. (47.43,261.27) .. controls (47.43,262.68) and (48.58,263.83) .. (50,263.83) .. controls (51.42,263.83) and (52.57,262.68) .. (52.57,261.27) -- cycle ;
\draw  [fill={rgb, 255:red, 0; green, 0; blue, 0 }  ,fill opacity=1 ] (142.83,261.27) .. controls (142.83,259.85) and (141.68,258.7) .. (140.27,258.7) .. controls (138.85,258.7) and (137.7,259.85) .. (137.7,261.27) .. controls (137.7,262.68) and (138.85,263.83) .. (140.27,263.83) .. controls (141.68,263.83) and (142.83,262.68) .. (142.83,261.27) -- cycle ;
\draw    (50,171) -- (140.27,261.27) ;
\draw    (50,261.27) -- (140.27,171) ;

\draw (38,151.67) node [anchor=north west][inner sep=0.75pt]   [align=left] {$x$};
\draw (37.33,259.33) node [anchor=north west][inner sep=0.75pt]   [align=left] {$y$};
\draw (142.27,264.27) node [anchor=north west][inner sep=0.75pt]   [align=left] {$z$};
\draw (142,151) node [anchor=north west][inner sep=0.75pt]   [align=left] {$u$};

\end{tikzpicture}

    \caption{ A complete graph with 4 vertices $K_4$}
    \label{example_5}
\end{figure}

{\it Example 5.} $G=(V,E,\mathbf{w}_0)$ is complete graph with 4 vertices $K_4$ in Figure \ref{example_5}, $V=\{x,y,z,u\}$, $E=\{xz,xy,xu,zy,zu,yu\}$, $\mathbf{w}_0=(w_0,w_0,w_0,w_0,w_0,w_0)$. For weights $\mathbf{w}=(w,w,w,w,w,w)$,
we have for $\alpha\in[0,1]$,
$$\mu_x^\alpha(v)=\left\{\begin{array}{lll}
\alpha,&v=x\\
\frac{1-\alpha}{3},&v=y,\,z,\,u,
\end{array}\right.\mu_y^\alpha(v)=\left\{\begin{array}{lll}
\alpha,&v=y\\
\frac{1-\alpha}{3},&v=x,\,z,\,u.
\end{array}\right.$$
Moreover, we calculate
$$W(\mu_x^\alpha,\mu_y^\alpha)-d(x,y)=\left\{\begin{array}{lll}
\frac{-2-4\alpha}{3}w,&\alpha\in[0,1/4]\\
\frac{4(\alpha-1)}{3}w,&\alpha\in[1/4,1].
\end{array}\right.$$
Setting $\mu(v,\cdot)=\mu_v(\cdot)$ for all $v\in V$, we have that the evolution problem (\ref{evolution}) has a unique global solution
$\mathbf{w}(t)=(w(t),w(t),w(t),w(t))$, where
$$w(t)=\left\{\begin{array}{lll}
w_0e^{-\frac{2+4\alpha}{3}t},&\alpha\in[0,1/4]\\
w_0e^{\frac{4(\alpha-1)}{3}t},&\alpha\in[1/4,1].
\end{array}\right.$$

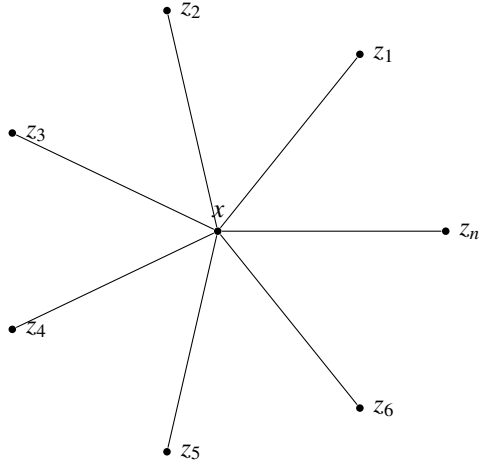
\begin{figure}[H]
    \centering
    \begin{tikzpicture}
    \def\n{7} 
    \def\radius{3} 

    \node[circle, fill=black, inner sep=1pt, label=above:$x$] (center) at (0,0) {};

    \foreach \i in {1,2,3,4,5,6} {
        \node[circle, fill=black, inner sep=1pt, label=right:$z_{\i}$] (z\i) at ({\i*360/\n}: \radius) {};
        \draw (center) -- (z\i);
    }

    \node[circle, fill=black, inner sep=1pt, label=right:$z_{n}$] (z7) at ({7*360/\n}: \radius) {};
    \draw (center) -- (z7);

\end{tikzpicture}

    \caption{A star with $n+1$ vertices}
    \label{example_6}
\end{figure}

{\it Example 6.} Let $G=(V,E,\mathbf{w}_0)$ be a star in Figure \ref{example_6}, $V=\{x,z_1,z_2,\cdots,z_n\}$, $E=\{xz_1,xz_2,\cdots,xz_n\}$,
$\mathbf{w}_0=(w_0,w_0,\cdots,w_0)$. For weights $\mathbf{w}=(w,w,\cdots,w)$,
we have for $\alpha\in[0,1]$,
$$\mu_x^\alpha(v)=\left\{\begin{array}{lll}
\alpha,&v=x\\
\frac{1-\alpha}{n},&v\in V\setminus\{x\},
\end{array}\right.\,\,\mu_{y_i}^\alpha(v)=\left\{\begin{array}{lll}
\alpha,&v=z_i\\
{1-\alpha},&v=x\\
0,& v\in V\setminus\{x,z_i\}.
\end{array}\right.$$
 We also have
$$W(\mu_x^\alpha,\mu_{z_i}^\alpha)=\left\{\begin{array}{lll}
(1-2\alpha)w,&\alpha\in[0,1/(n+1)]\\
\frac{n+2\alpha-2}{n}w,&\alpha\in[1/(n+1),1],
\end{array}\right.$$
and thus
$$W(\mu_x^\alpha,\mu_{z_i}^\alpha)-d(x,z_i)=\left\{\begin{array}{lll}
-2\alpha w,&\alpha\in[0,1/(n+1)]\\
\frac{2(\alpha-1)}{n}w,&\alpha\in[1/(n+1),1],\end{array}\right.$$
$i=1,2,\cdots,n$.
Setting $\mu(v,\cdot)=\mu_v(\cdot)$ for all $v\in V$, we derive that the flow (\ref{evolution}) has a unique global solution
$\mathbf{w}(t)=(w(t),w(t),w(t),w(t))$, where
$$w(t)=\left\{\begin{array}{lll}
w_0e^{-2\alpha t},&\alpha\in[0,1/n]\\
w_0e^{\frac{2(\alpha-1)}{n}t},&\alpha\in[1/n,1].
\end{array}\right.$$

 \section{Experimental results}\label{experiment}
In this section, we first discretize systems (\ref{evolution}) and (\ref{evolution-norm}), and then apply the corresponding algorithms to community detection. Next, we use the usual criteria and real-world datasets to evaluate the accuracy of our method in community detection. Then, we examine how some key parameters influence the accuracy of algorithm. Finally, our approach is compared with existing methods.
\subsection{Discretization and algorithm design}
If we choose $\mu(x,z)$ to be $\alpha$-lazy one-step random walk $\mu_x^\alpha(z)$ (\ref{alpha-lazy}), then the system (\ref{evolution})
is a modification of (\ref{Olloivier-flow}), namely
\begin{equation}\label{Ollivier}
\left\{\begin{array}{lll}
w_{e_i}^\prime(t)=-\kappa_{e_i}^\alpha(t) \rho_{e_i}(t)\\[1.2ex]
e_i=x_iy_i, \, \kappa_{e_i}^\alpha(t)=1-\frac{W(\mu_{x_i}^\alpha,\,\mu_{y_i}^\alpha)}{\rho_{e_i}(t)}\\[1.2ex]
w_{e_i}(0)=w_{0,i},\,\,
i=1,2,\cdots,m.
\end{array}\right.
\end{equation}
While the system (\ref{evolution-norm}) becomes a quasi-normalized version of (\ref{Ollivier}), which reads as
\begin{equation}\label{Ollivier-norm}
\left\{\begin{array}{lll}
w_{e_i}^\prime(t)=-\kappa_{e_i}^\alpha(t) \rho_{e_i}(t)+\frac{\sum_{j=1}^m \kappa_{e_j}(t)\rho_{e_j}(t)}{\sum_{j=1}^mw_{e_j}(t)}\rho_{e_i}(t)\\[1.2ex]
e_i=x_iy_i, \, \kappa_{e_i}^\alpha(t)=1-\frac{W(\mu_{x_i}^\alpha,\,\mu_{y_i}^\alpha)}{\rho_{e_i}(t)}\\[1.2ex]
w_{e_i}(0)=w_{0,i},\,\,
i=1,2,\cdots,m.
\end{array}\right.
\end{equation}

If we choose $\mu(x,z)$ as an $\alpha$-lazy two-step random walk given by (\ref{2-lazy}), then we have special cases of
(\ref{evolution}) and (\ref{evolution-norm}) as follows:
\begin{equation}\label{evolution-2step}
\left\{\begin{array}{lll}
w_{e_i}^\prime(t)=W(\mu_{2,x_i}^\alpha,\mu_{2,y_i}^\alpha)-d(x_i,y_i)\\[1.2ex]
e_i=x_iy_i\in E\\[1.2ex]
w_{e_i}(0)=w_{0,i},\,\,
i=1,2,\cdots,m,
\end{array}\right.
\end{equation}

\begin{equation}\label{evolution-norm-2step}
\left\{\begin{array}{lll}
w_{e_i}^\prime(t)=W\left(\mu_{2,x_i}^\alpha,\mu_{2,y_i}^\alpha\right)-d(x_i,y_i)-\frac{\sum_{j=1}^m
\left(W(\mu_{2,x_j}^\alpha,\mu_{2,y_j}^\alpha)-d(x_j,y_j)\right)}
{\sum_{j=1}^mw_{e_j}} d(x_i,y_i)\\[1.2ex]
e_i=x_iy_i\in E\\[1.2ex]
w_{e_i}(0)=w_{0,i},\,\,
i=1,2,\cdots,m.
\end{array}\right.
\end{equation}

A discrete version of (\ref{Ollivier}) can be written as
$$w_{e_i}(t+s)-w_{e_i}(t)=-s\kappa_{e_i}^\alpha(t) \rho_{e_i}(t),$$
where $s>0$ is a step size of discretization. Since other discrete versions of (\ref{Ollivier-norm})-(\ref{evolution-norm-2step}) are very similar, we omitted them here. To balance the computer's calculation accuracy and error, we set step size $s=0.01$.

 The pseudo-code in Algorithm \ref{one_algorithm} below demonstrates the iterations process for the discrete versions of (\ref{Ollivier}) and (\ref{Ollivier-norm}), and the surgery in the final iteration.

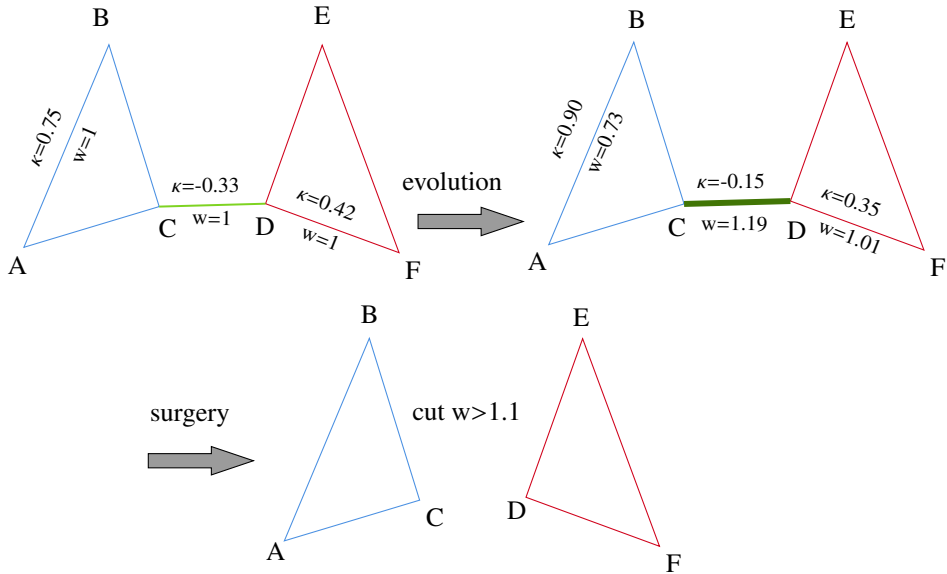
\begin{figure}[t]
    \centering
\begin{tikzpicture}[x=0.75pt,y=0.75pt,yscale=-1,xscale=1]

\draw  [color={rgb, 255:red, 74; green, 144; blue, 226 }  ,draw opacity=1 ] (342.54,103.52) -- (275.01,123.99) -- (317.44,22.11) -- cycle ;
\draw  [color={rgb, 255:red, 208; green, 2; blue, 27 }  ,draw opacity=1 ] (395.64,102.06) -- (423.91,22.3) -- (462.27,126.72) -- cycle ;
\draw [color={rgb, 255:red, 65; green, 117; blue, 5 }  ,draw opacity=1 ][line width=2.25]    (342.54,103.53) -- (395.63,102.06) ;
\draw  [color={rgb, 255:red, 0; green, 0; blue, 0 }  ,draw opacity=1 ][fill={rgb, 255:red, 155; green, 155; blue, 155 }  ,fill opacity=1 ] (209.76,108.79) -- (241.25,108.79) -- (241.25,105.22) -- (262.24,112.36) -- (241.25,119.49) -- (241.25,115.92) -- (209.76,115.92) -- cycle ;
\draw  [color={rgb, 255:red, 74; green, 144; blue, 226 }  ,draw opacity=1 ] (210.54,252.52) -- (143.01,272.99) -- (185.44,171.11) -- cycle ;
\draw  [color={rgb, 255:red, 208; green, 2; blue, 27 }  ,draw opacity=1 ] (263.64,251.06) -- (291.91,171.3) -- (330.27,275.72) -- cycle ;
\draw  [color={rgb, 255:red, 0; green, 0; blue, 0 }  ,draw opacity=1 ][fill={rgb, 255:red, 155; green, 155; blue, 155 }  ,fill opacity=1 ] (75.1,229.12) -- (106.58,229.12) -- (106.58,225.55) -- (127.57,232.69) -- (106.58,239.83) -- (106.58,236.26) -- (75.1,236.26) -- cycle ;
\draw  [color={rgb, 255:red, 74; green, 144; blue, 226 }  ,draw opacity=1 ] (80.54,104.85) -- (13.01,125.33) -- (55.44,23.44) -- cycle ;
\draw  [color={rgb, 255:red, 208; green, 2; blue, 27 }  ,draw opacity=1 ] (133.64,103.4) -- (161.91,23.64) -- (200.27,128.05) -- cycle ;
\draw [color={rgb, 255:red, 126; green, 211; blue, 33 }  ,draw opacity=1 ][line width=0.75]    (80.54,104.86) -- (133.63,103.4) ;

\draw (33.91,81.83) node [anchor=north west][inner sep=0.75pt]  [rotate=-293] [align=left] {{\footnotesize w=1}};
\draw (95.81,105.39) node [anchor=north west][inner sep=0.75pt]  [rotate=-0.52] [align=left] {{\footnotesize w=1}};
\draw (152.6,111.8) node [anchor=north west][inner sep=0.75pt]  [rotate=-21.7] [align=left] {{\footnotesize w=1}};
\draw (85.56,89.33) node [anchor=north west][inner sep=0.75pt]   [align=left] {{\footnotesize \( {\kappa } \)=-0.33}};
\draw (12.29,82.87) node [anchor=north west][inner sep=0.75pt]  [rotate=-293] [align=left] {{\footnotesize \( {\kappa } \)=0.75}};
\draw (150.79,92.47) node [anchor=north west][inner sep=0.75pt]  [rotate=-21.71] [align=left] {{\footnotesize \( {\kappa } \)=0.42}};
\draw (291.12,86.78) node [anchor=north west][inner sep=0.75pt]  [rotate=-293.9] [align=left] {{\footnotesize w=0.73}};
\draw (349.91,108.26) node [anchor=north west][inner sep=0.75pt]  [rotate=-0.52] [align=left] {{\footnotesize w=1.19}};
\draw (410.67,108.77) node [anchor=north west][inner sep=0.75pt]  [rotate=-21.67] [align=left] {{\footnotesize w=1.01}};
\draw (347.83,86.95) node [anchor=north west][inner sep=0.75pt]   [align=left] {{\footnotesize \( {\kappa } \)=-0.15}};
\draw (272.45,78.75) node [anchor=north west][inner sep=0.75pt]  [rotate=-293] [align=left] {{\footnotesize \( {\kappa } \)=0.90}};
\draw (412.56,91.42) node [anchor=north west][inner sep=0.75pt]  [rotate=-20.22] [align=left] {{\footnotesize \( {\kappa } \)=0.35}};
\draw (200.78,86.12) node [anchor=north west][inner sep=0.75pt]   [align=left] {evolution};
\draw (205.64,202.05) node [anchor=north west][inner sep=0.75pt]   [align=left] {cut w$>$1.1};
\draw (74.78,205.12) node [anchor=north west][inner sep=0.75pt]   [align=left] {surgery};
\draw (3.68,128.33) node [anchor=north west][inner sep=0.75pt]   [align=left] {A};
\draw (132.67,272.66) node [anchor=north west][inner sep=0.75pt]   [align=left] {A};
\draw (263.33,126.33) node [anchor=north west][inner sep=0.75pt]   [align=left] {A};
\draw (46,3.67) node [anchor=north west][inner sep=0.75pt]   [align=left] {B};
\draw (79.17,110) node [anchor=north west][inner sep=0.75pt]   [align=left] {C};
\draw (127.33,109) node [anchor=north west][inner sep=0.75pt]   [align=left] {D};
\draw (156.67,2.67) node [anchor=north west][inner sep=0.75pt]   [align=left] {E};
\draw (202.27,131.05) node [anchor=north west][inner sep=0.75pt]   [align=left] {F};
\draw (212.54,255.52) node [anchor=north west][inner sep=0.75pt]   [align=left] {C};
\draw (333.33,109.33) node [anchor=north west][inner sep=0.75pt]   [align=left] {C};
\draw (252.97,252.06) node [anchor=north west][inner sep=0.75pt]   [align=left] {D};
\draw (392.67,111) node [anchor=north west][inner sep=0.75pt]   [align=left] {D};
\draw (313.33,5) node [anchor=north west][inner sep=0.75pt]   [align=left] {B};
\draw (180,153.33) node [anchor=north west][inner sep=0.75pt]   [align=left] {B};
\draw (286,154.67) node [anchor=north west][inner sep=0.75pt]   [align=left] {E};
\draw (418.67,5) node [anchor=north west][inner sep=0.75pt]   [align=left] {E};
\draw (332.27,278.72) node [anchor=north west][inner sep=0.75pt]   [align=left] {F};
\draw (464.27,129.72) node [anchor=north west][inner sep=0.75pt]   [align=left] {F};

\end{tikzpicture}

    \caption{Demonstration of evolution of weights on a graph.}
    \label{fig:one_evolution}
\end{figure}

\begin{algorithm}[H]

        \caption{Community detection using $\alpha$-lazy one-step random walk evolution }
        \label{one_algorithm}
        \KwIn{an undirected connected finite network \( G = \left( {V,E}\right)  \) ; parameter \(\alpha\) ; maximum iteration \( T \) ; step size \( s \) .}
        \KwOut{community detection  results of $G$}

       \For{ \( i = 1,\cdots ,T \)}
       {\, calculate the distance $d(x,y)$ between any two points $x$ and $y$;\\
        compute the Ricci curvature \( {\kappa }_{e}^{\alpha} \) for all edge;\\update all edge weights simultaneously through the corresponding flow process;\\
        }
        \For{cutoff= $w_{max},\cdots, w_{min}$}
        {
        \For{ \( e \in  E \)}{
            \If{$w_e>{\it cutoff}$}{
      remove the edge $e$;
      }
            }
        calculate the accuracy of community detection;
        }
\end{algorithm}
Now we demonstrate a simple application of Algorithm 1. Let $G$ be as in Figure \ref{fig:one_evolution}, $w$ be the weight of edge and $\kappa$ be Ollivier's Ricci curvature given by (\ref{Ollivier-norm}) with $\alpha=0.5$. Red and blue represent the edges within the two communities, and green edges are the edges connecting the two communities. The initial weight on each edge is assumed to be 1. The weights and curvatures of edges not shown can be derived from the symmetry. The first graph is the initial calculated curvature and weight. The second graph is the result after thirty iterations by (\ref{Ollivier-norm}). The bottom graph is the community detection result after the surgery.

Similar to Algorithm \ref{one_algorithm}, the discrete versions of (\ref{evolution-2step}) and (\ref{evolution-norm-2step}) can be applied to design Algorithm \ref{two_algorithm}.
\begin{algorithm}[H]{}

        \caption{Community detection using $\alpha$-lazy two-step random walk evolution }
        \label{two_algorithm}
        \KwIn{an undirected connected finite network \( G = \left( {V,E}\right)  \) ; parameter \(\alpha\) ; maximum iteration \( T \) ; step size \( s \) .}
        \KwOut{community detection  results of $G$}

       \For{ \( i = 1,\cdots ,T \)}{
        \, calculate the distance $d(x,y)$ between any two points $x$ and $y$;\\
        compute the Wasserstein distance  $W(\mu_{2,x_i}^\alpha,\mu_{2,y_i}^\alpha)$ for all edge;\\update all edge weights simultaneously through the corresponding flow process;\\
        }
        \For{cutoff= $w_{max},\cdots, w_{min}$}
        {
        \For{ \( e \in  E \)}{
            \If{$w_e>{\it cutoff}$}{
      remove the edge $e$;
      }
            }
        calculate the accuracy of community detection;
        }
\end{algorithm}

Our algorithms are primarily complex due to the tasks of finding the shortest path in the graph and solving a linear programming problem. The time complexities for these tasks are \(O(|E|\log|V|)\) and \(O(|E|D^3)\) respectively. Here, \(D\) represents the average degree of the network, and $|E|$ and $|V|$ represent the number of edges and vertices in the network, respectively. Despite the sparse connectivity of the network, where \(|D| \ll |E|\), \(O(|E|D^3)\) often surpasses \(O(|E|\log|V|)\) in most scenarios. Consequently, the computational complexity of our approach is \(O(|E|D^3)\).

\subsection{Real-world datasets and criteria}

 For the real-world datasets, we select three distinct scale community graphs. The three real-world benchmark networks used in our experiments are drawn from the datasets introduced in \cite{Ni-Lin}. The results obtained in the present setting are comparable with those obtained when all edge weights are set to one; for brevity, we omit the details here and leave them to interested readers. We shall use three different metrics to assess the precision of community detection in real-world datasets. We choose the adjusted Rand index (ARI) \cite{Hubert-Arabie} and normalized mutual information (NMI) \cite{Danon-Guilera-Duch} as the criteria for evaluating the quality of clustering accuracy when compared to the ground truth. Furthermore, we choose modularity \cite{Clauset-Newman-Moore, M.Newman} to measure the robustness of the community structure of a given graph without relying on ground-truth clustering. Basic information for real-world networks is listed in Table \ref{tab:1}.
\begin{table}[t]
\tbl{Summary of real-world network characteristics\label{tab:1}}
{\begin{tabular}{@{}ccccccc@{}}
\toprule
Networks & Vertices & Edges & \#Class & AvgDeg & Density  &Diameter\\
\colrule
Karate   & 34      & 78    & 2                     & 4.59   & 0.139  &5 \\
Football & 115     & 613   & 12                    & 10.66  & 0.094   &4\\
Facebook & 775     & 14006 & 18                    & 36.15  & 0.047    &	8\\
\botrule
\end{tabular}}
\end{table}

The Karate dataset, referenced as \cite{Zachary}, was collected from the members of a university karate club by Wayne Zachary in the 1970s. This network consists of 34 members and 78 edges, where the vertices represent individual members and the edges denote the connections between them. This data set is generally used to find the two groups of people into which the karate club fission after a conflict between two faculties.

In the football dataset, referenced as \cite{Girvan M}, the focus is on the NCAA football Division I games schedule for the 2000 season. It includes 115 teams (vertices) and 613 matches (edges). Each vertex represents a specific team, and the edges indicate the matches played between these teams. The dataset identifies 12 ground-truth communities or conferences. Since teams tend to play more frequently within their own conference, this network clearly exhibits a community structure.

After analyzing the NCAA Football League Network, we expanded our analysis to three larger-scale networks. These networks present additional challenges, particularly in terms of computational efficiency. The Facebook network dataset, known as \cite{Jure L}, is derived from the Stanford Network Analysis Project and captures interaction networks within an online social networking platform. The benchmark community ground truth in this dataset is meticulously organized based on well-defined themes such as interests and affiliations.

Next, we describe the three criteria for evaluation employed in our assessment process. To be more specific,  we let
\(\{ U_1, U_2, \ldots, U_I \} \) and \(\{ W_1, W_2, \ldots, W_J \} \) be two partitions of the set \( S \) of $n$ vertices (nodes). Denote \( m_{ij} =|U_i\cap W_j|\) the number of vertices in \( U_i\cap W_j \), while \( c_i \) and \( d_j \) represent the numbers of vertices in \( U_i \) and \( W_j \), respectively. All these quantities are listed in Table \ref{tab:2}.
\begin{table}[t]
\tbl{Contingency table\label{tab:2}}
{\begin{tabular}{c|cccc|c}
$U \backslash W$         & $W_1$    & $W_2$    & $\cdots$ & $W_J$    & \text{sums} \\ \hline
$U_1$                    & $m_{11}$ & $m_{12}$ & $\cdots$ & $m_{1J}$ & $c_1$                    \\
$U_2$                    & $m_{21}$ & $m_{22}$ & $\cdots$ & $m_{2J}$ & $c_2$                    \\
$\vdots$                 & $\vdots$ & $\vdots$ & $\ddots$ & $\vdots$ & $\vdots$                 \\
$U_I$                    & $m_{I1}$ & $m_{I2}$ & $\cdots$ & $m_{IJ}$ & $c_I$                    \\ \hline
\text{sums} & $d_1$    & $d_2$    & $\cdots$ & $d_J$    &
\end{tabular}}
\end{table}

 Then the explicit expressions of the above mentioned three criteria are written below.

\begin{itemize}
    \item \textbf{Adjusted Rand Index}
    \[
    \text{ARI} = \frac{\sum_{i=1}^{I}\sum_{j=1}^{J} \binom{m_{ij}}{2} - \left( \sum_{i=1}^{I} \binom{c_i}{2} \sum_{j=1}^{J} \binom{d_j}{2} \big/ \binom{n}{2} \right)}{\frac{1}{2} \left( \sum_{i=1}^{I} \binom{c_i}{2} +
    \sum_{j=1}^{J} \binom{d_j}{2} \right) - \left( \sum_{i=1}^{I} \binom{c_i}{2} \sum_{j=1}^{J} \binom{d_j}{2} \big/ \binom{n}{2} \right)},
    \]
    where \( \binom{n}{2} \) is the number of different pairs from \( n \) vertices, while symbols \( \binom{m_{ij}}{2} \), \( \binom{c_i}{2} \) and \( \binom{d_j}{2} \) have the same meaning. The Adjusted Rand Index (ARI) is a statistical measure used in data clustering analysis. It quantifies the similarity between two partitions of a dataset by comparing the assignments of data points to clusters. The ARI value ranges from -1 to 1, where a value of 1 indicates a perfect match between the partitions and a value close to 0 indicates a random assignment of data points to clusters.

    \item \textbf{Normalized Mutual Information}
    \[
    \text{NMI} = \frac{-2 \sum_{i=1}^I \sum_{j=1}^J m_{ij} \log \left( \frac{m_{ij} n}{c_i d_j} \right)}{\sum_{i=1}^I c_i \log \left( \frac{c_i}{n} \right) + \sum_{j=1}^J d_j \log \left( \frac{d_j}{n} \right)}.
    \]
The Normalized Mutual Information (NMI) is a metric commonly utilized to evaluate the effectiveness of community detection algorithms. It facilitates the comparison of two clusters or communities by producing a value within the range of 0 to 1. A higher NMI value signifies a stronger resemblance between the two partitions or communities.

    \item \textbf{Modularity}
    \[
    Q = \sum_{k=1}^M \left( \frac{C_k}{|E|} - \beta \left( \frac{D_k}{2|E|} \right)^2 \right),
    \]
 where \( M \) represents the number of communities, \( C_k \) is the number of connections within the $k$th community, \( D_k \) is the total degree of vertices in the $k$th community, and \( \beta \) is a resolution parameter, with a default value of \(1\). The parameter \( Q \) spans from \(-0.5\) to \(1\), with values approaching \(1\) signifying a more robust community structure and superior partition quality. A positive value indicates an excess of edges within groups compared to what would be expected by chance. Modularity assesses how edges connect within specific groups compared to random link distribution among all nodes.
\end{itemize}

The ARI is a pair-counting measure that compares the number of element pairs correctly grouped together or separated in two partitions. It adjusts for the likelihood that some agreement might occur by chance, thereby offering a more reliable comparison between partitions than raw accuracy scores. Its sensitivity to cluster size differences makes it a particularly useful tool in cases where clusters vary significantly in size.

NMI, in contrast, originates from information theory, where it captures the shared information content between two partitions. By normalizing this shared information with respect to the entropy of each partition, NMI ensures consistency even as the number of communities changes, providing robustness to fluctuations in community structure.

Modularity Q offers an edge-centric perspective. It measures how well a network is partitioned into communities by comparing the observed fraction of intra-community edges with what would be expected in a random network. A higher Q value implies a more distinct community structure, independent of ground truth labels. This makes modularity particularly valuable in exploratory analyses where community structure is unknown.

To ensure clarity and coherence in the discussion, the following terminology will be used throughout:

\begin{itemize}
    \item \textbf{one\_evol} designates the modified Ollivier’s Ricci flow (\ref{Ollivier}), incorporating an $\alpha$-lazy one-step random walk.
    \item The quasi-normalized form of this $\alpha$-lazy one-step random walk (\ref{Ollivier-norm} )will be referred to as \textbf{qn1\_evol}.
    \item When the system (\ref{evolution-2step}) utilizes an $\alpha$-lazy two-step random walk, it will be labeled \textbf{two\_evol}.
    \item The quasi-normalized variant of two-step system (\ref{evolution-norm-2step}) will be called \textbf{qn2\_evol}.
\end{itemize}

\subsection{Analysis of key parameters}
In Algorithms \ref{one_algorithm} and \ref{two_algorithm}, parameters \( \alpha \) and the maximum iteration \( T \) significantly influence the final results. In this section, we shall examine the impact of these two parameters on the experimental outcomes, detailed in Sections \ref{var_alpha} and \ref{var_iteration}, respectively.
\subsubsection{Impact of \texorpdfstring{$\alpha$}{} on experimental results}\label{var_alpha}
\begin{figure}[H]
    \centering
    \includegraphics[width=0.8\linewidth]{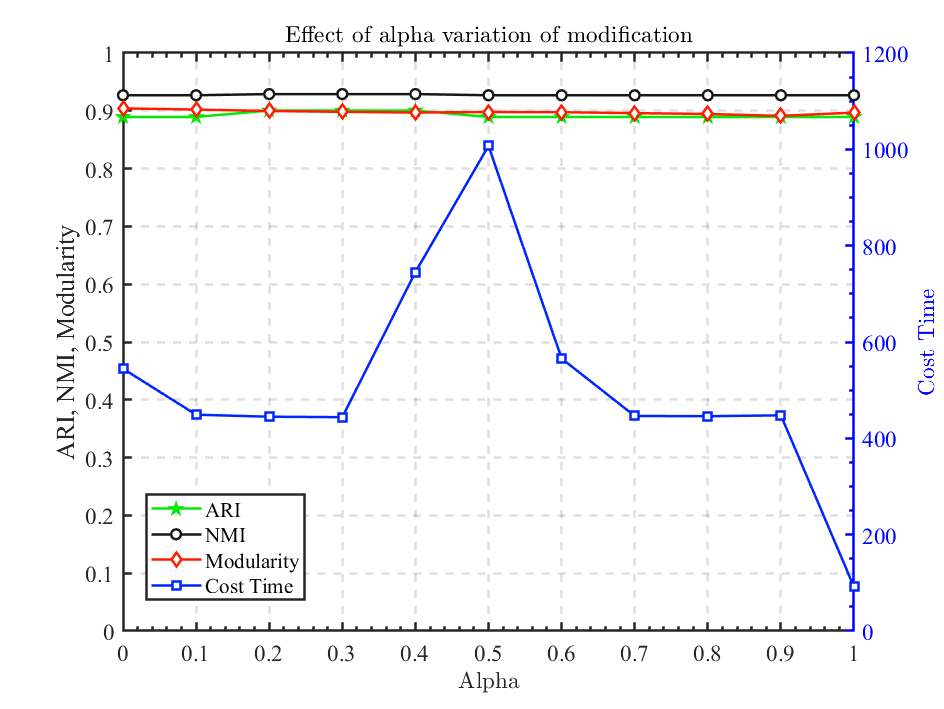}
    \caption{\textbf{two\_evol} with different alpha setting on football network}
    \label{alpha_1}
\end{figure}

\begin{figure}[htbp]
    \centering
    \includegraphics[width=0.8\linewidth]{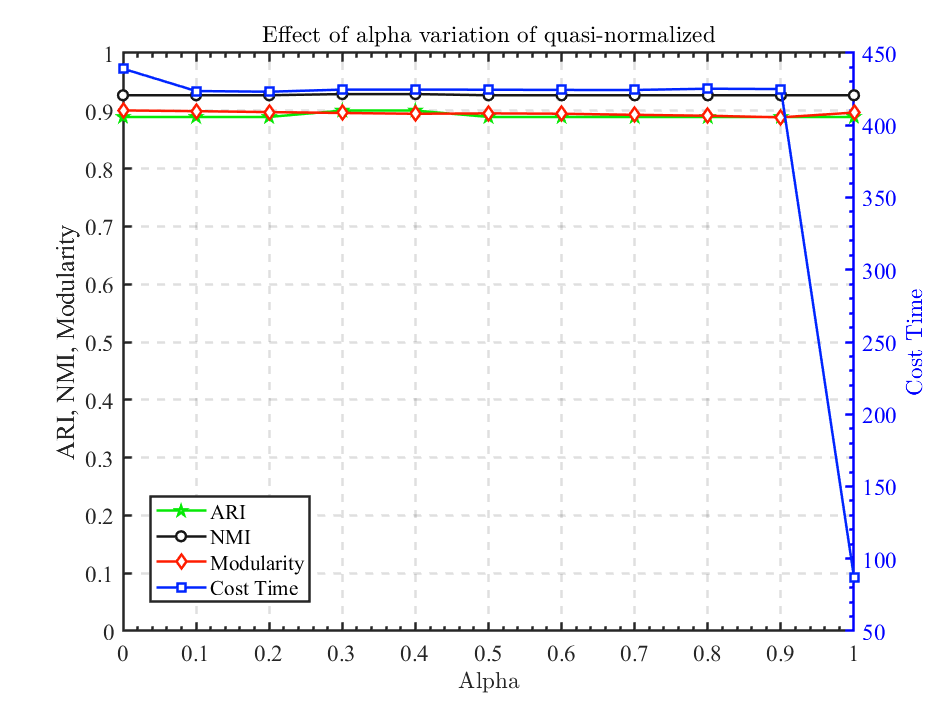}
    \caption{\textbf{qn2\_evol} with different alpha setting on football network}
    \label{alpha_2}
\end{figure}

Note that all of the aforementioned systems rely on the selection of $\alpha$. Extensive experiments conducted by C. C. Ni et al \cite{Ni-Lin} have shown that when applying community detection based on (\ref{Olloivier-flow}) with an $\alpha$-lazy one-step random walk $\mu_x^\alpha(z)$, $\alpha = 0.5$ yields the best results. Next, we shall study the impact of the alpha value on the $\alpha$-lazy two-step random walk system (\ref{evolution-2step}) and (\ref{evolution-norm-2step}). The value range of $\alpha$ is limited to $[0,1)$. The experimental results are listed in Figures \ref{alpha_1} and \ref{alpha_2}.

In Figure \ref{alpha_1}, the graph shows that ARI remains stable around 0.9 across different $\alpha$ values, indicating that the clustering quality produced by the {two\_evol} algorithm does not significantly change with variations in $\alpha$. This suggests a high degree of consistency in clustering performance. Similar to ARI, the NMI values are almost constant at around 0.9, indicating that the clustering similarity remains high regardless of the $\alpha$ value, and thus the algorithm produces highly stable clustering results. The Modularity line shows minimal fluctuation, maintaining a value close to 0.9 across all $\alpha$ values. This indicates that the quality of community detection within the network remains steady even as $\alpha$ changes, and the modularity is robust to parameter variations. The cost time shows a distinct variation from $\alpha$. At $\alpha = 0.5$, the cost time peaks sharply, reaching around 1000 ms. However, as $\alpha$ approaches $1$, the cost time rapidly decreases to nearly 100 ms.

In Figure \ref{alpha_2}, ARI remains very stable, consistently near 0.9, indicating minimal impact of $\alpha$ on the clustering quality in the {qn2\_evol} algorithm. NMI also stays steady at approximately 0.9, similar to the first graph, further reinforcing the stability of the algorithm with respect to clustering similarity. Modularity shows little variation, remaining close to 0.9 across different $\alpha$ values, suggesting that the quality of community detection in {qn2\_evol} is largely unaffected by changes in $\alpha$. Compared to {two\_evol}, the cost time for {qn2\_evol} is notably lower, hovering around 400 ms for most $\alpha$ values. Similar to {two\_evol}, the computation time significantly decreases to nearly zero as $\alpha$ approaches 1, although the overall computational cost for {qn2\_evol} remains much lower across the entire range of $\alpha$ values. This suggests that {qn2\_evol} is generally more efficient than {two\_evol}, especially in terms of computational cost, making it potentially more suitable for applications where efficiency is a priority.

The ARI, NMI, and modularity curves in both graphs remain consistently high and stable (around 0.9), indicating that both algorithms maintain strong clustering performance across all $\alpha$ values.
This stability suggests that the clustering quality and the effectiveness of network partitioning are not highly sensitive to changes in the $\alpha$ parameter. Therefore, the choice of $\alpha$ does not substantially affect the algorithms' ability to produce accurate and consistent clustering.

Cost time varies significantly with $\alpha$, showing a sharp peak for {two\_evol} around $\alpha$ = 0.5, while {qn2\_evol} displays lower and more stable computational costs. For both algorithms, the computational cost drops sharply at $\alpha$ approaches 1, implying that the algorithms become much more efficient at this value. However, {qn2\_evol} consistently exhibits lower computational complexity than {two\_evol}, suggesting that it may be a more efficient algorithm overall.

The sharp drop in computation time at $\alpha$ approaches 1 suggests that this particular value of $\alpha$ could be leveraging some form of internal optimization in both algorithms, drastically reducing the time complexity. When $\alpha$ approaches 1, the expression in (\ref{2-lazy}) becomes relatively straightforward. This phenomenon provides a guideline for selecting $\alpha$ values that optimize computational efficiency without sacrificing clustering performance.\\

 In the subsequent sections, we stick to $\alpha=0.5$ as our experiment parameter setting to maintain consistency in algorithms for $\alpha$-lazy one-step random walk and two-step random walk.

\subsubsection{Effect of iterations on experimental results}\label{var_iteration}
The impact of iterations on experimental results is illustrated in Figures \ref{fig:one_iteration} through \ref{fig:qn2_iteration} below. Our analysis indicates that the application of normalization had a negligible impact on the final results of community detection, regardless of whether it was implemented. In other words, Figures \ref{fig:one_iteration} and \ref{fig:qn1_iteration} display the same pattern, while Figures \ref{fig:two_iteration} and \ref{fig:qn2_iteration} are nearly identical. Based on this observation, we focus comparison on the qn1\_evol and qn2\_evol algorithms.

\begin{figure}[htbp!]
    \centering
    \begin{subfigure}[t]{0.32\textwidth}
        \centering
        \includegraphics[width=\linewidth]{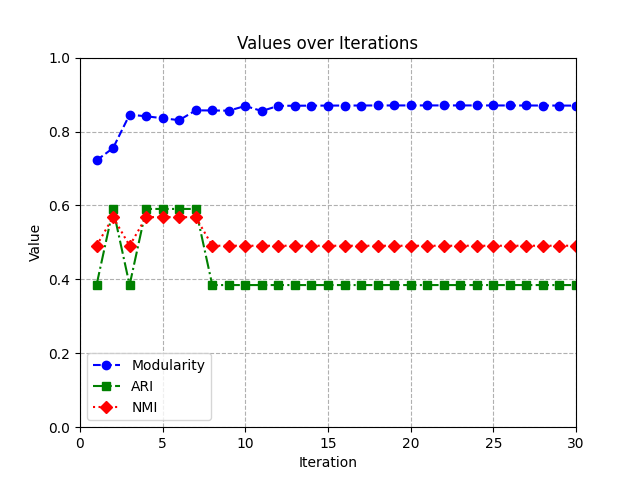}
        \caption{karate}
        \label{fig:one_iteration_karate}
    \end{subfigure}
    \hfill
    \begin{subfigure}[t]{0.32\textwidth}
        \centering
        \includegraphics[width=\linewidth]{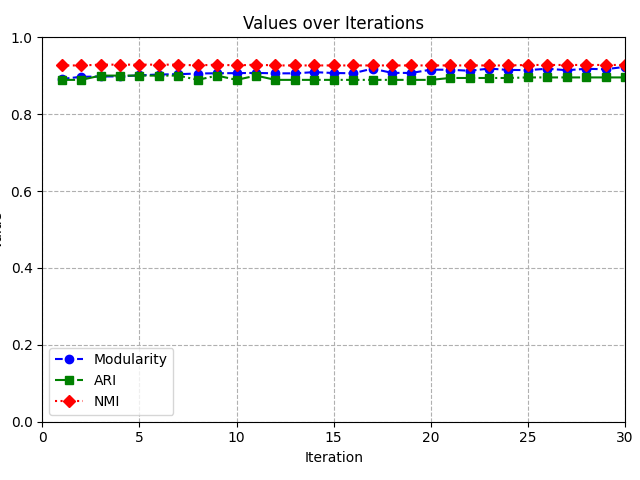}
        \caption{football}
        \label{fig:one_iteration_football}
    \end{subfigure}
    \hfill
    \begin{subfigure}[t]{0.32\textwidth}
        \centering
        \includegraphics[width=\linewidth]{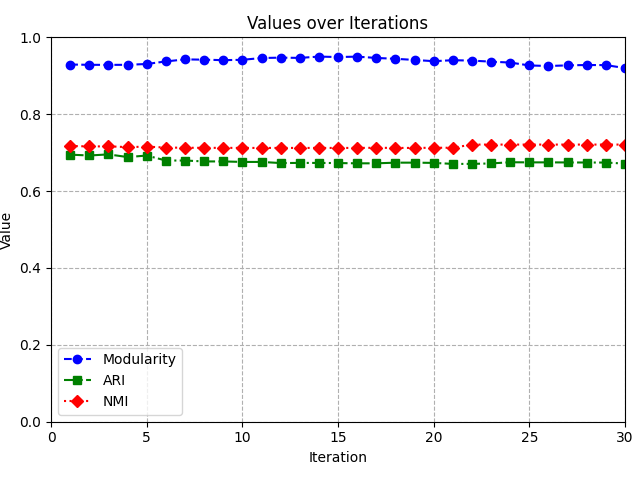}
        \caption{facebook}
        \label{fig:one_iteration_facebook}
    \end{subfigure}
    \caption{\textbf{one\_evol} over iterations}
    \label{fig:one_iteration}
\end{figure}

\begin{figure}[htbp!]
    \centering
    \begin{subfigure}[t]{0.32\textwidth}
        \centering
        \includegraphics[width=\linewidth]{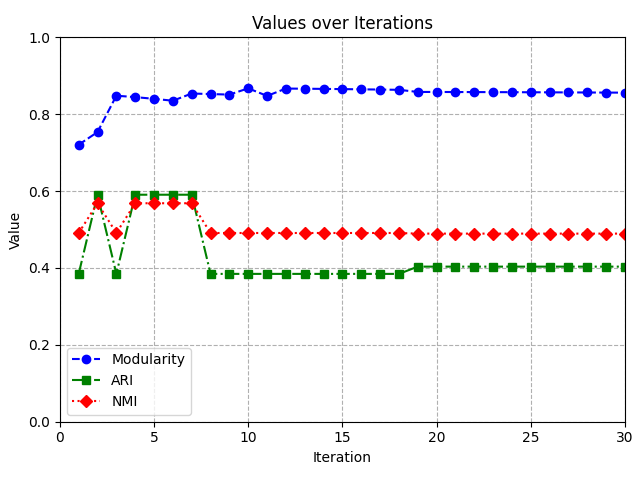}
        \caption{karate}
        \label{fig:qn1_iteration_karate}
    \end{subfigure}
    \hfill
    \begin{subfigure}[t]{0.32\textwidth}
        \centering
        \includegraphics[width=\linewidth]{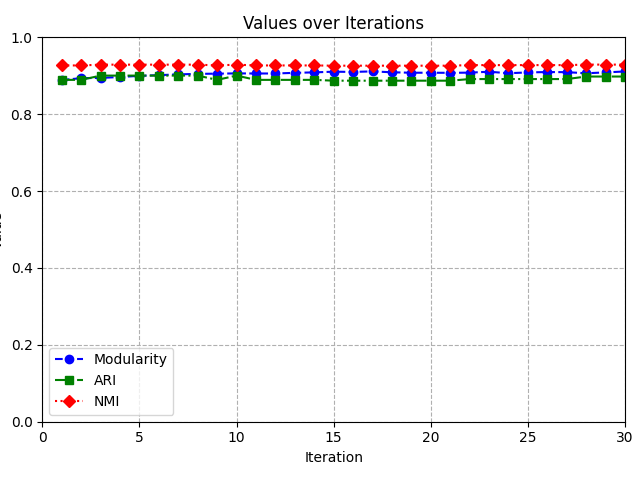}
        \caption{football}
        \label{fig:qn1_iteration_football}
    \end{subfigure}
    \hfill
    \begin{subfigure}[t]{0.32\textwidth}
        \centering
        \includegraphics[width=\linewidth]{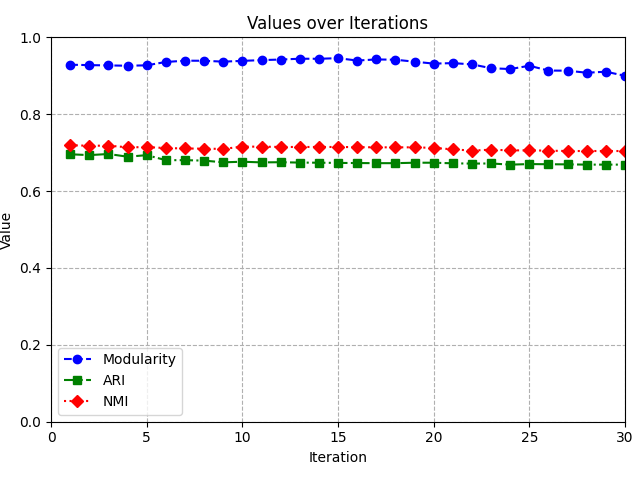}
        \caption{facebook}
        \label{fig:qn1_iteration_facebook}
    \end{subfigure}
    \caption{\textbf{qn1\_evol} over iterations}
    \label{fig:qn1_iteration}
\end{figure}

In Figure \ref{fig:qn1_iteration} (a), modularity starts below 0.8, increases over the first few iterations, and stabilizes around 0.85. ARI and NMI have some fluctuations in the initial few iterations, with ARI settling around 0.48 and NMI stabilizing at around 0.4.
(b) indicates all three metrics (modularity, ARI, NMI) are stable throughout the iterations, starting at higher values (modularity $\approx$ 0.9, ARI $\approx$ 0.85, NMI $\approx$ 0.8) and maintaining these levels across all iterations.
While (c) shows that modularity starts slightly above 0.8 and stabilizes around 0.85, while ARI and NMI remain stable throughout (ARI $\approx$ 0.85 and NMI $\approx$ 0.75).

\begin{figure}[htbp!]
    \centering
    \begin{subfigure}[t]{0.32\textwidth}
        \centering
        \includegraphics[width=\linewidth]{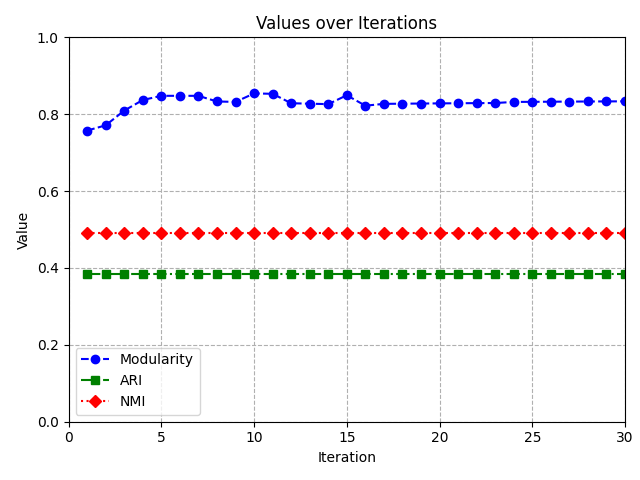}
        \caption{karate}
        \label{fig:two_iteration_karate}
    \end{subfigure}
    \hfill
    \begin{subfigure}[t]{0.32\textwidth}
        \centering
        \includegraphics[width=\linewidth]{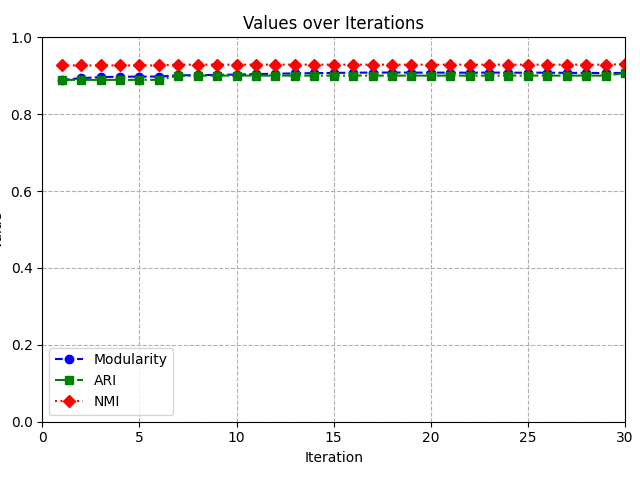}
        \caption{football}
        \label{fig:two_iteration_football}
    \end{subfigure}
    \hfill
    \begin{subfigure}[t]{0.32\textwidth}
        \centering
        \includegraphics[width=\linewidth]{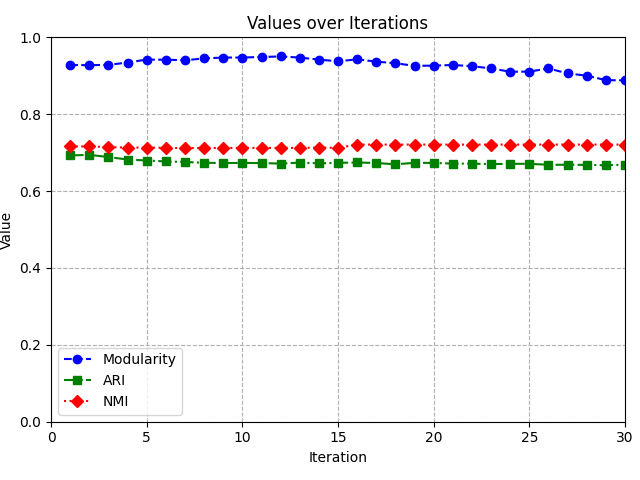}
        \caption{facebook}
        \label{fig:two_iteration_facebook}
    \end{subfigure}
    \caption{\textbf{two\_evol} over iterations}
    \label{fig:two_iteration}
\end{figure}

\begin{figure}[htbp!]
    \centering
    \begin{subfigure}[t]{0.32\textwidth}
        \centering
        \includegraphics[width=\linewidth]{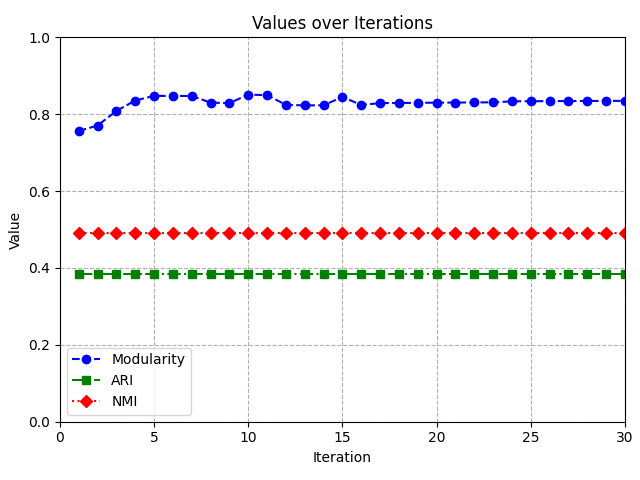}
        \caption{karate}
        \label{fig:qn2_iteration_karate}
    \end{subfigure}
    \hfill
    \begin{subfigure}[t]{0.32\textwidth}
        \centering
        \includegraphics[width=\linewidth]{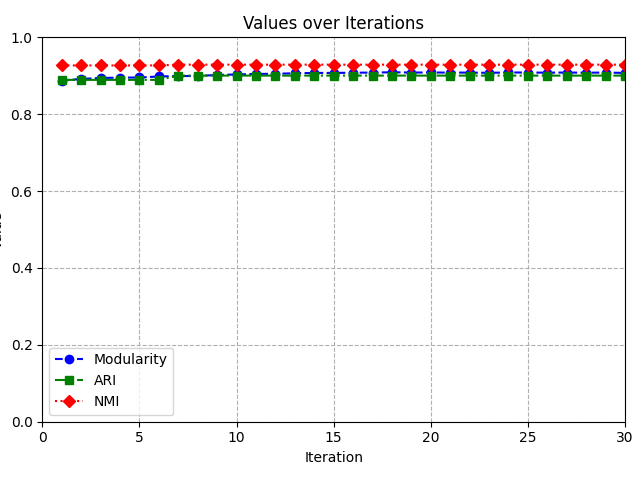}
        \caption{football}
        \label{fig:qn2_iteration_football}
    \end{subfigure}
    \hfill
    \begin{subfigure}[t]{0.32\textwidth}
        \centering
        \includegraphics[width=\linewidth]{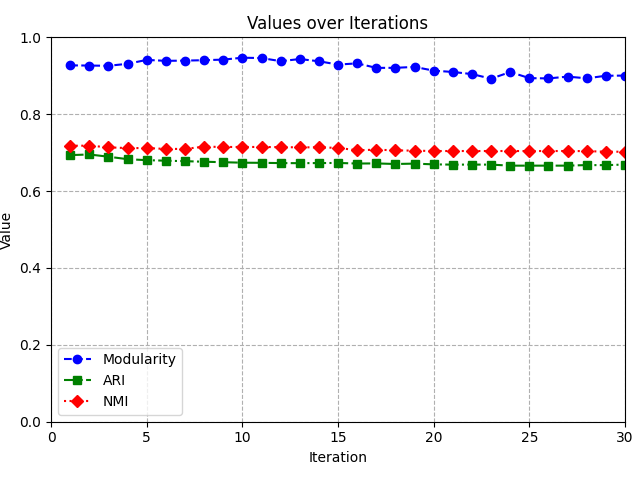}
        \caption{facebook}
        \label{fig:qn2_iteration_facebook}
    \end{subfigure}
    \caption{\textbf{qn2\_evol} over iterations}
    \label{fig:qn2_iteration}
\end{figure}

In Figure \ref{fig:qn2_iteration} (a), modularity starts around 0.8, increases and stabilizes around 0.85 after some oscillation. ARI and NMI exhibit minor fluctuations at the beginning, with ARI settling around 0.6 and NMI around 0.45.
(b) shows all metrics are stable, with Modularity $\approx$ 0.9, ARI $\approx$  0.85, and NMI $\approx$  0.8 from the beginning until the end of the iterations.
(c) says modularity starts at about 0.9 and gradually decreases to stabilize around 0.85. ARI and NMI remain stable throughout, with ARI $\approx$  0.85 and NMI $\approx$  0.75.

When evaluating qn1\_evol on the small karate network in Figure \ref{fig:qn1_iteration} (a), the modularity score exhibited rapid growth in the initial iterations, reaching a stable value of approximately 0.85 after five iterations. This suggests that the algorithm can efficiently optimize modularity, demonstrating a capacity for quick convergence. Such convergence implies that the algorithm is effective in detecting community structures within the network. However, the adjusted Rand index (ARI) metric displays noticeable volatility during the early stages, with an initial increase followed by a sharp decline, eventually stabilizing at a lower value near 0.4. This indicates that the algorithm undergoes substantial changes early on, leading to an initial improvement in performance but later succumbing to overfitting, which results in diminished accuracy in the long term. In practical scenarios, this issue could be mitigated by fine-tuning the number of iterations. Similarly, the normalized mutual information (NMI) measure follows a comparable pattern to ARI, experiencing minor fluctuations at the start but generally stabilizing around 0.5.

For the qn2\_evol algorithm on the same karate network in Figure \ref{fig:qn2_iteration} (a), the modularity indicator shows a slower increase compared to qn1\_evol, but it eventually stabilizes at a similar final value of 0.85. This more gradual rise can be attributed to the algorithm's conservative optimization strategy, likely incorporating regularization to prevent excessive fluctuations caused by rapid adjustments. While the slower growth may seem less efficient, the steadier performance implies greater robustness in optimizing modularity. The ARI score for this algorithm remains consistently around 0.4, with little variation, suggesting a stable, though less adaptable, approach. Similarly, the NMI score exhibits minimal change from the outset, hovering around 0.5, indicating limitations in the algorithm's ability to maintain cluster consistency in smaller networks.

In the evaluation of medium-sized networks, such as the football network depicted in Figures \ref{fig:qn1_iteration} (b) and  \ref{fig:qn2_iteration} (b), as well as large networks, exemplified by Facebook in Figures \ref{fig:qn1_iteration} (c) and \ref{fig:qn2_iteration} (c), both qn1\_evol and qn2\_evol demonstrate comparable performance across the three primary metrics. In both cases, the algorithms converge on optimal results after relatively few iterations, underscoring their efficiency in community detection. Despite minor fluctuations in certain indicators as the number of iterations increases, both algorithms demonstrate strong overall consistency, highlighting their stability and robustness in handling networks of varying sizes.

In summary, tests on networks of different sizes show that the qn1\_evol and qn2\_evol algorithms work similarly for community detection. Both algorithms reach their results quickly and require fewer iterations, which makes them efficient and stable. Although there are slight changes in the indicators during the process, they maintain overall consistency, demonstrating their reliability. The qn1\_evol algorithm optimizes modularity faster in small networks but can overfit. In contrast, the qn2\_evol algorithm takes a more careful approach, providing better stability. These findings suggest that both algorithms perform well in community detection for various network sizes, especially in medium to large networks.

\subsection{\texorpdfstring{$\alpha$}{}-lazy one-step  random walk}\label{section_5.1}
In this subsection, the experimental results are shown in various charts comprising three parts:
(a) histograms of Ricci curvatures and edge weights before discrete flow;
(b) histograms of Ricci curvatures and edge weights after discrete flow;
(c) evaluation metrics after surgery. The results of the $\alpha$-lazy one-step random walk on three real-world datasets are presented from Figures \ref{fig:one_karate} to  \ref{fig:one_facebook} below.

\begin{figure}[htbp!]
    \centering
    \begin{subfigure}[t]{0.32\textwidth}
        \centering
        \includegraphics[width=\linewidth]{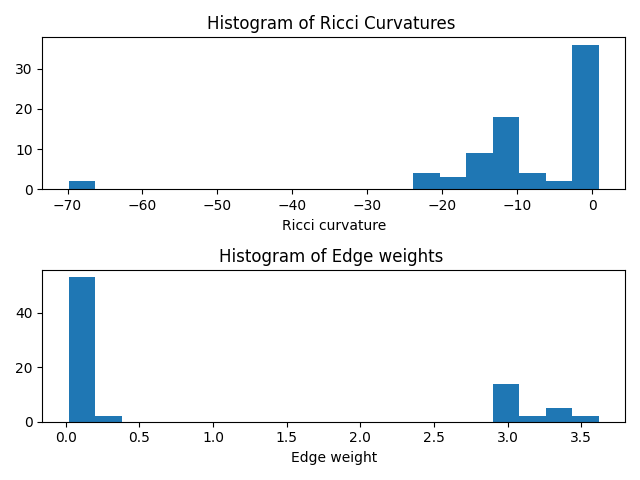}
        \caption{}
        \label{fig:one_karate_1}
    \end{subfigure}
    \hfill
    \begin{subfigure}[t]{0.32\textwidth}
        \centering
        \includegraphics[width=\linewidth]{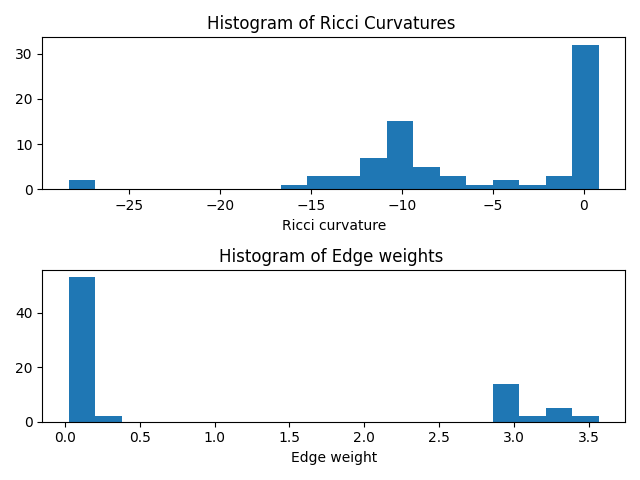}
        \caption{}
        \label{fig:one_karate_2}
    \end{subfigure}
    \hfill
    \begin{subfigure}[t]{0.32\textwidth}
        \centering
        \includegraphics[width=\linewidth]{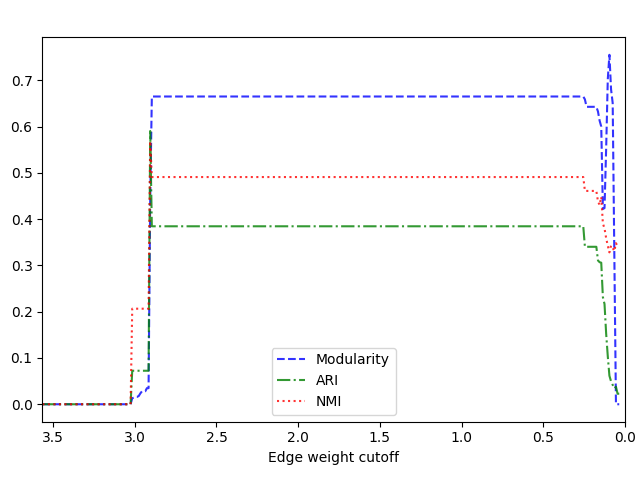}
        \caption{}
        \label{fig:one_karate_3}
    \end{subfigure}
    \caption{one\_evol on karate}
    \label{fig:one_karate}
\end{figure}

In Figure \ref{fig:one_karate} (a), the initial curvature distribution across edges spans from approximately -70 to 0, with most values clustering near 0 and a few showing highly negative curvatures, around -70. These significantly negative curvatures may highlight regions within the graph where structural imbalance is prominent. Additionally, most edge weights are near zero, while a smaller subset approaches values of 3 to 3.5, signaling varying connectivity strengths across the graph.

After five iterations, Figure \ref{fig:one_karate} (b) illustrates a markedly narrower curvature distribution, now ranging between -25 and 0. The Ricci flow process effectively mitigates extreme negative curvatures, bringing more edges close to zero curvature, which suggests an improved structural balance. The edge weight distribution, though still skewed toward low values, demonstrates a slight shift toward higher weights, indicating that iterative weight adjustments may contribute to stabilizing graph connectivity.

Figure \ref{fig:one_karate} (c) tracks the behavior of various metrics as the edge weight threshold decreases. Initially, three metrics show consistency at higher thresholds, rising to an optimal point before dropping off as the cutoff approaches zero. This suggests that excessive removal of low-weight edges can disrupt community structure, emphasizing the importance of balanced pruning in preserving the graph’s integrity.

\begin{figure}[htbp!]
    \centering
    \begin{subfigure}[t]{0.32\textwidth}
        \centering
        \includegraphics[width=\linewidth]{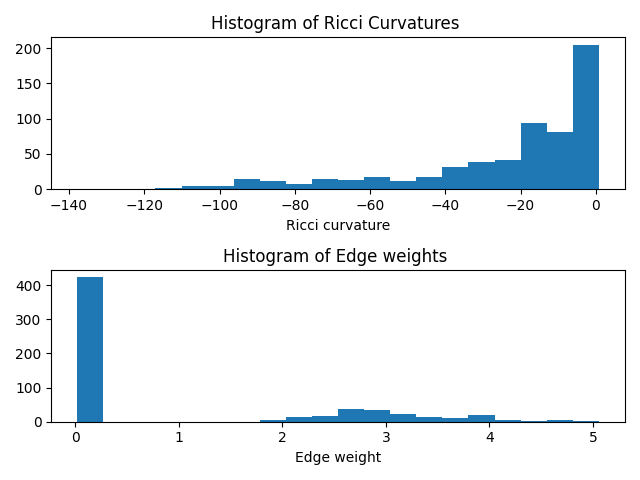}
        \caption{}
        \label{fig:one_football_1}
    \end{subfigure}
    \hfill
    \begin{subfigure}[t]{0.32\textwidth}
        \centering
        \includegraphics[width=\linewidth]{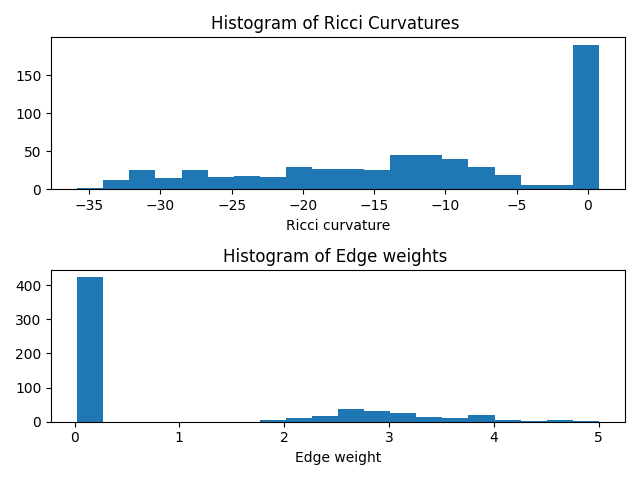}
        \caption{}
        \label{fig:one_football_2}
    \end{subfigure}
    \hfill
    \begin{subfigure}[t]{0.32\textwidth}
        \centering
        \includegraphics[width=\linewidth]{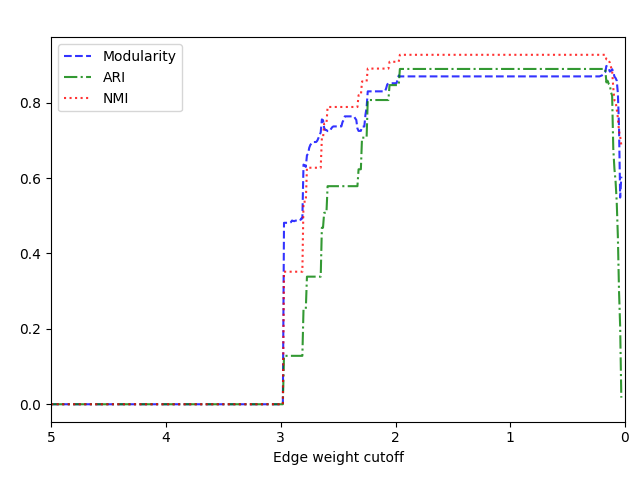}
        \caption{}
        \label{fig:one_football_3}
    \end{subfigure}
    \caption{one\_evol on football}
    \label{fig:one_football}
\end{figure}

In Figure \ref{fig:one_football} (a), initial Ricci curvature values range widely from approximately -140 to 0, with most edges exhibiting curvatures near 0 and a smaller portion reaching significantly negative values around -140. This distribution suggests areas of substantial geometric imbalance within the graph, as these highly negative curvatures indicate regions of structural tension or weak connectivity. The distribution of edge weights is also skewed, with most weights close to zero, and a few reaching approximately 5, which implies a sparse, weakly connected network before applying the discrete flow.

After the Ricci flow, as illustrated in Figure \ref{fig:one_football} (b), the curvature values become more compressed, spanning a reduced range from -35 to 0. This reduction in extreme negative values and more even curvature distribution signal that the  flow has effectively smoothed the graph’s geometry, diminishing pronounced imbalances. Though edge weights remain predominantly low, there is a minor shift towards higher weights, indicating a subtle increase in connectivity strength and enhanced structural balance across the graph.

In Figure \ref{fig:one_football} (c), modularity sharply rises as edge weight cutoffs decrease from 5 to around 3, then stabilizes at an elevated level, suggesting that removing weaker edges sharpens community clarity. Beyond a weight threshold of about 3, further pruning yields diminishing modularity gains. A similar trend is seen in the adjusted Rand index (ARI), which increases with edge removal, aligning with the effects seen in normalized mutual information (NMI). NMI increases as weak edges are eliminated but stabilizes near 0.9, before declining when edge removal is too aggressive. This suggests that excessive pruning may compromise the integrity of community structure, emphasizing the need for balanced edge retention to preserve overall connectivity.

\begin{figure}[htbp]
    \centering
    \begin{subfigure}[t]{0.32\textwidth}
        \centering
        \includegraphics[width=\linewidth]{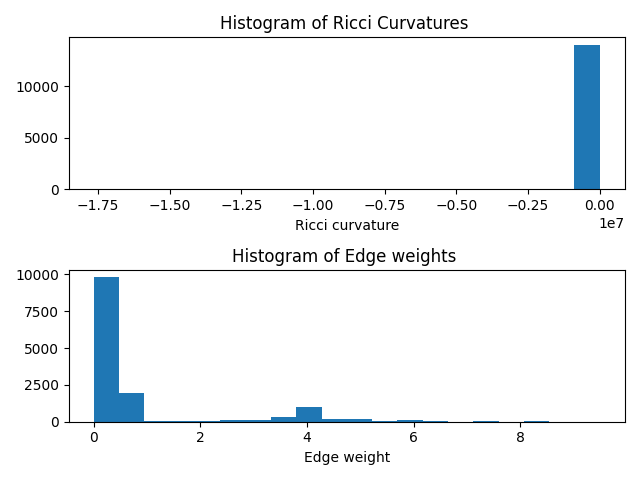}
        \caption{}
        \label{fig:one_facebook_1}
    \end{subfigure}
    \hfill
    \begin{subfigure}[t]{0.32\textwidth}
        \centering
        \includegraphics[width=\linewidth]{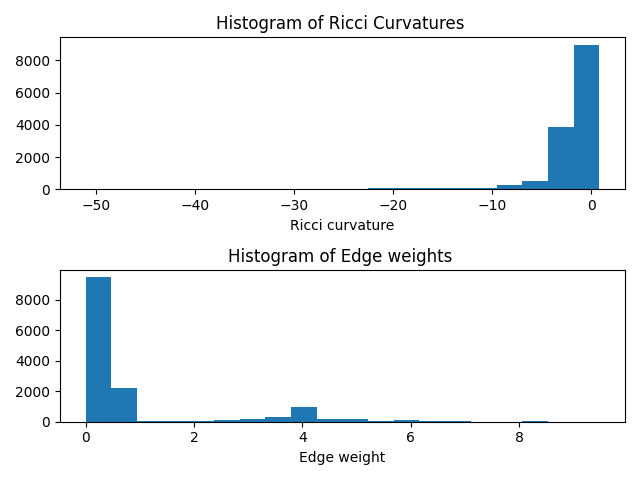}
        \caption{}
        \label{fig:one_facebook_2}
    \end{subfigure}
    \hfill
    \begin{subfigure}[t]{0.32\textwidth}
        \centering
        \includegraphics[width=\linewidth]{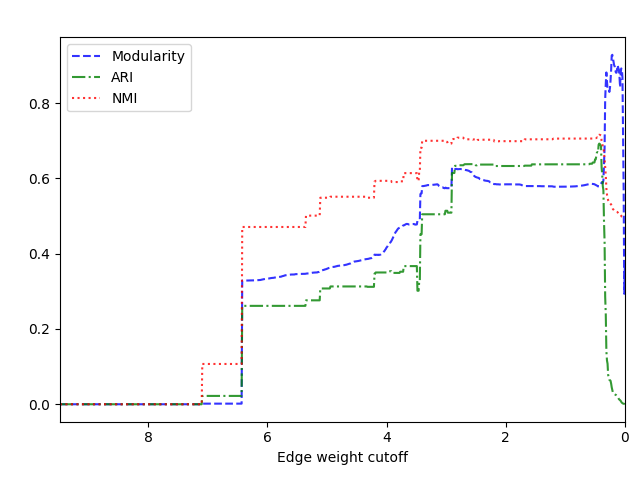}
        \caption{}
        \label{fig:one_facebook_3}
    \end{subfigure}
    \caption{one\_evol on facebook}
    \label{fig:one_facebook}
\end{figure}

In Figure \ref{fig:one_facebook} (a), the initial Ricci curvature distribution is notably narrow, with most values clustering near zero, suggesting a largely uniform structural balance in the Facebook graph. A limited number of edges have slightly negative curvatures, yet the overall variation is minimal before the flow application. The edge weight distribution is skewed toward lower values, with the majority of weights near zero and only a few higher-weight edges reaching up to 8 or 9. This distribution implies that, although the graph has a few strong connections, most edges represent weak connections prior to the flow.

After Ricci flow, as depicted in Figure \ref{fig:one_facebook} (b), the range of Ricci curvatures expands considerably, now extending from -50 to 0. This broader distribution introduces more negative curvature values, suggesting that the Ricci flow enhances geometric diversity and potentially strengthens structural balance within the graph. Edge weights remain largely consistent with their initial spread; low-weight edges are still prevalent, while a small portion maintains higher weights. This pattern implies that, with limited iterations, the Ricci flow mainly adjusts curvature while leaving the edge weight distribution largely unchanged, reflecting a refined structural equilibrium.

Figure \ref{fig:one_facebook} (c) illustrates that modularity increases as weak edges are pruned, signifying enhanced community clarity. However, further pruning beyond certain thresholds yields diminishing returns, with modularity eventually stabilizing. The adjusted Rand index (ARI) also rises as weaker edges are removed but drops sharply at very low cutoffs, indicating a sensitivity of community structure to aggressive edge removal. Normalized mutual information (NMI) follows a similar trend, rising with edge pruning and stabilizing over a range of cutoffs, suggesting that moderate pruning strengthens community integrity. However, excessive removal ultimately weakens the structure, underscoring the need for balanced edge retention to maintain cohesive community structures.

\subsection{\texorpdfstring{$\alpha$}{}-lazy two-step  random walk}\label{section_5.2}
\begin{figure}[htbp!]
    \centering
    \begin{subfigure}[t]{0.32\textwidth}
        \centering
        \includegraphics[width=\linewidth]{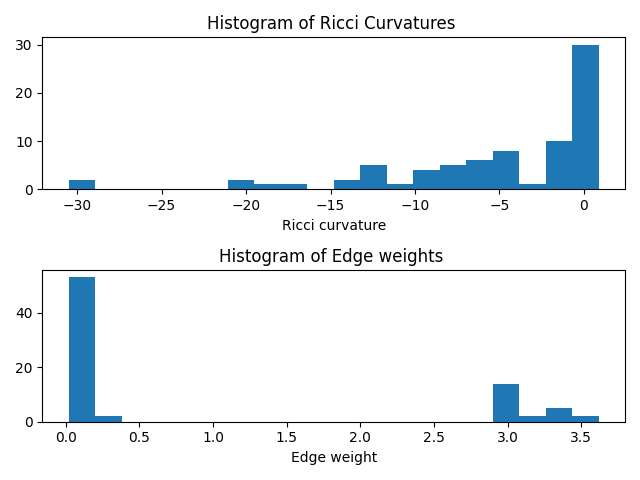}
        \caption{}
        \label{fig:two_karate_1}
    \end{subfigure}
    \hfill
    \begin{subfigure}[t]{0.32\textwidth}
        \centering
        \includegraphics[width=\linewidth]{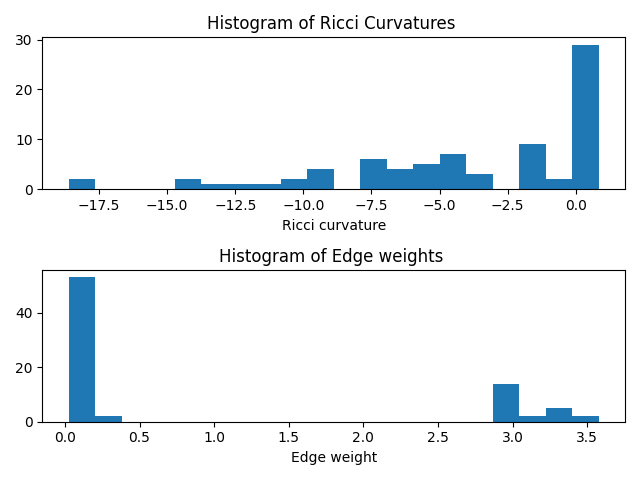}
        \caption{}
        \label{fig:two_karate_2}
    \end{subfigure}
    \hfill
    \begin{subfigure}[t]{0.32\textwidth}
        \centering
        \includegraphics[width=\linewidth]{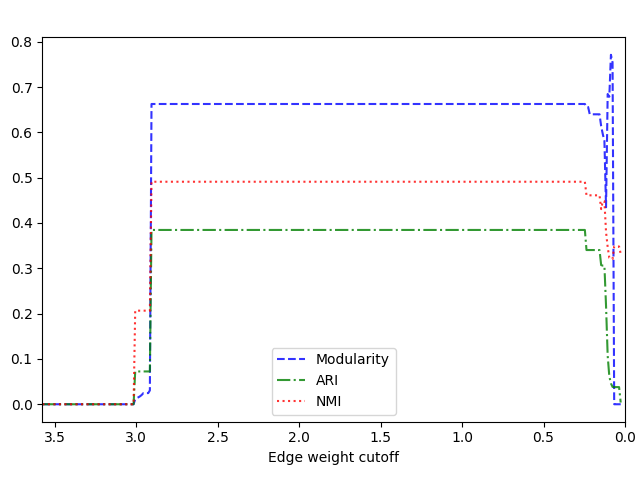}
        \caption{}
        \label{fig:two_karate_3}
    \end{subfigure}
    \caption{two\_evol on karate}
    \label{fig:two_karate}
\end{figure}

\begin{figure}[htbp!]
    \centering
    \begin{subfigure}[t]{0.32\textwidth}
        \centering
        \includegraphics[width=\linewidth]{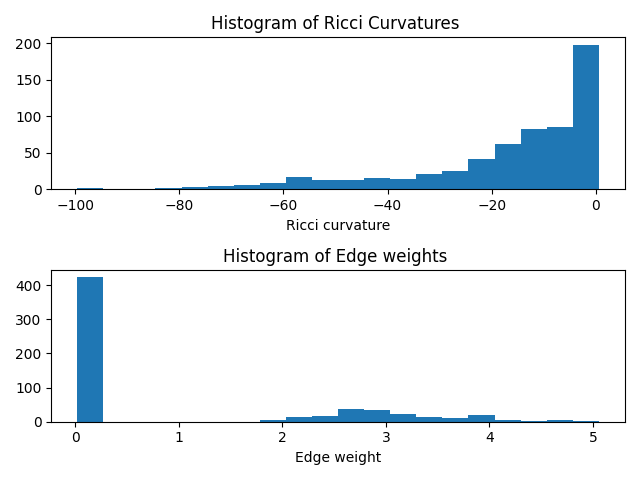}
        \caption{}
        \label{fig:two_football_1}
    \end{subfigure}
    \hfill
    \begin{subfigure}[t]{0.32\textwidth}
        \centering
        \includegraphics[width=\linewidth]{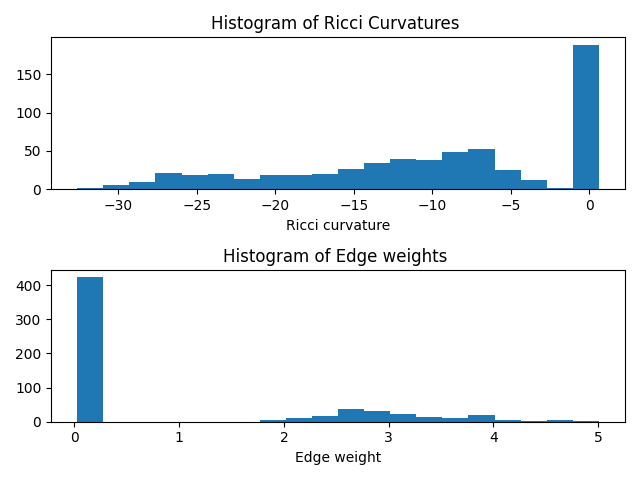}
        \caption{}
        \label{fig:two_football_2}
    \end{subfigure}
    \hfill
    \begin{subfigure}[t]{0.32\textwidth}
        \centering
        \includegraphics[width=\linewidth]{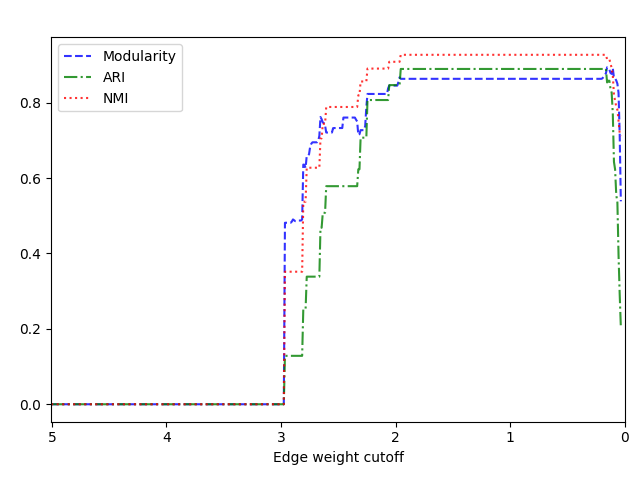}
        \caption{}
        \label{fig:two_football_3}
    \end{subfigure}
    \caption{two\_evol on football}
    \label{fig:two_football}
\end{figure}

\begin{figure}[htbp!]
    \centering
    \begin{subfigure}[t]{0.32\textwidth}
        \centering
        \includegraphics[width=\linewidth]{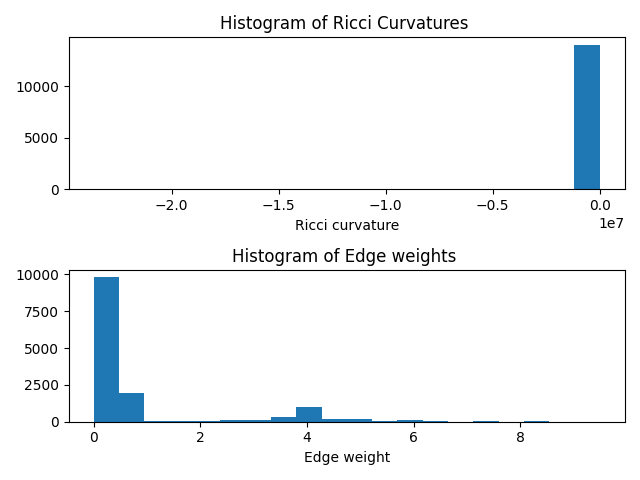}
        \caption{}
        \label{fig:two_facebook_1}
    \end{subfigure}
    \hfill
    \begin{subfigure}[t]{0.32\textwidth}
        \centering
        \includegraphics[width=\linewidth]{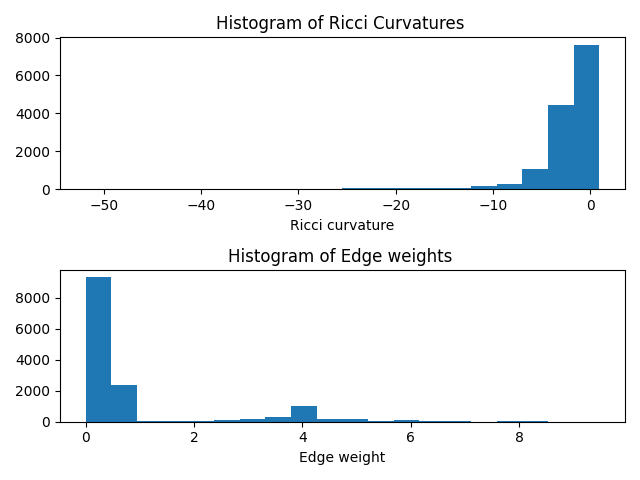}
        \caption{}
        \label{fig:two_facebook_2}
    \end{subfigure}
    \hfill
    \begin{subfigure}[t]{0.32\textwidth}
        \centering
        \includegraphics[width=\linewidth]{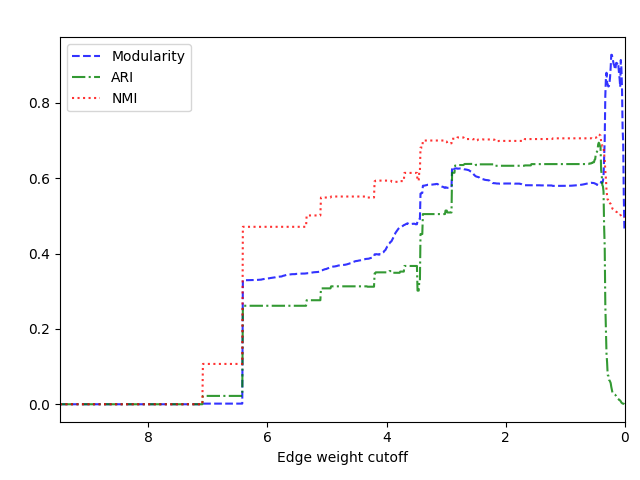}
        \caption{}
        \label{fig:two_facebook_3}
    \end{subfigure}
    \caption{two\_evol on facebook}
    \label{fig:two_facebook}
\end{figure}

Based on equations (\ref{evolution-2step}), (\ref{evolution-norm-2step}), and Algorithm \ref{two_algorithm}, we only focus on the distance $d(x_i,y_i)$ and Wasserstein distance $W(\mu_{2,x_i}^\alpha,\mu_{2,y_i}^\alpha)$, where $x_iy_i=e_i$ belongs to $E$. To maintain consistency with previous discussions, throughout this subsection, we continue to use the symbol of Ricci  curvature $$\kappa_{e_i}^\alpha=1-\frac{W(\mu_{2,x_i}^\alpha,\mu_{2,y_i}^\alpha)}{{d(x_i,y_i)}}.$$  It is important to note that introducing Ricci curvature is not essential in this context. The results of the $\alpha$-lazy two-step random walk on three real-world datasets are presented from Figures \ref{fig:two_karate} to \ref{fig:two_facebook}. Since the experimental results for the three datasets exhibit similar trends, we only discuss the experimental results of two\_evol on football. Figure \ref{fig:two_football} presents the histograms of Ricci curvatures and edge weights before and after applying the discrete flow, as well as the changes in modularity, ARI, and NMI with respect to the edge weight cutoff. Readers interested in further analysis will find similar patterns in Figures \ref{fig:two_karate} and \ref{fig:two_facebook}, suggesting a consistent effect of two\_evol across different networks.

Before applying the discrete Flow, both the distribution of Ricci curvatures and the edge weights exhibit notable imbalances in Figure \ref{fig:two_football} (a). The Ricci curvatures are broadly distributed, with most values concentrated in the negative range between -100 and 0, indicating significant negative curvature across the network. This suggests that the community structures in the initial network are sharply divided, with uneven geometric characteristics. Similarly, the edge weights show a concentration near zero, pointing to weak associations across much of the network. Only a few edges have higher weights, reflecting stronger connections in certain local areas, but overall, the network is characterized by sparse and weakly connected components.

After the application of Discrete Ricci Flow, the network structure undergoes a significant transformation. According to Figure \ref{fig:two_football} (b), the curvature distribution becomes more balanced, with values mostly ranging between -35 and 0, indicating a homogenization of the network’s geometric properties. This suggests that the flow has smoothed out the community boundaries, leading to a more cohesive and uniform structure. In parallel, the edge weight distribution remains centered around lower values, but there is an increase in the number of edges with higher weights. This shift indicates that Discrete Ricci Flow not only balances the network but also strengthens the associations within it, particularly in areas that previously had weak connections.

In terms of evaluation metrics in Figure \ref{fig:two_football} (c), modularity, Adjusted Rand Index (ARI), and Normalized Mutual Information (NMI) all display clear improvements following the application of Discrete Ricci Flow, as illustrated in the corresponding plots. As the edge weight cutoff decreases, modularity rapidly increases before stabilizing, indicating a more defined and pronounced community structure. Similarly, ARI shows a sharp increase before leveling off, reflecting improved clustering consistency, where the detected community structures align more closely with the true community divisions. Finally, the NMI follows a comparable trajectory, with its values rising and stabilizing at a high level, highlighting the enhanced information-theoretic similarity between the observed and actual community structures. These results collectively demonstrate that Discrete Ricci Flow not only optimizes the network's geometric structure but also improves its functional characteristics, as evidenced by the refined community detection metrics.

Initially, all these indicators are zero because only a few edges are removed, and there is minimal community structure. However, the indicators gradually increase, reach a peak, and then stabilize as the {\it cutoff} value increases. When the cutoff approaches $w_{min}$, most of the edges are deleted, leading to a rapid drop in these indicators, essentially reaching zero.

The application of discrete Ricci Flow leads to a significant improvement in the network’s geometric structure, as evidenced by the more uniform distribution of Ricci curvatures and the better-balanced edge weight distribution. The evaluation metrics demonstrate that community detection improves substantially as the edge weights are adjusted, with the metrics stabilizing at high values. This confirms the effectiveness of discrete Ricci Flow in enhancing both the performance and consistency of community detection algorithms.

\subsection{Comparison with other methods}\label{section_5.4}
 We shall compare our methods of community detection with three traditional machine learning methods, namely the Girvan Newman algorithm based on edge betweenness \cite{Girvan M}, the greedy modularity algorithm based on modularity maximization \cite{Clauset-Newman-Moore,Reichardt-Bornholdt}, and the label propagation algorithm based on stochastic methods \cite{Cordasco-Gargano}. We also use  another five different models based on Ricci curvature for comparison, including unnormalized discrete Ollivier's Ricci flow (DORF) \cite{Ni-Lin}, normalized discrete Ollivier's Ricci flow (NDORF), normalized discrete Lin-Lu-Yau's Ricci flow (NDSRF) \cite{Lai X}, and modifications of discrete Lin-Lu-Yau's Ricci flow (Rho; RhoN) \cite{M-Y1}. The results in Table \ref{tab:3} demonstrate the effectiveness of our evolutionary algorithms (the last four methods) in detecting communities in real-world scenarios.

\begin{table}[t]
\tbl{Performance of various algorithms on real-world networks\label{tab:3}}
{\begin{tabular}{@{}cccccccccc@{}}
\toprule
Methods\textbackslash{}Networks & \multicolumn{3}{c}{Karate club}               & \multicolumn{3}{c}{Football}                  & \multicolumn{3}{c}{Facebook}                  \\
                                                 & ARI           & NMI           & Q             & ARI           & NMI           & Q             & ARI           & NMI           & Q    \\         \colrule
Girvan Newman                                    & \textbf{0.77} & \textbf{0.73} & 0.48          & 0.14          & 0.36          & 0.50          & 0.03          & 0.16          & 0.01          \\
Greedy Modularity                                & 0.57          & 0.56          & 0.58          & 0.47          & 0.70          & 0.82          & 0.49          & 0.68          & 0.55          \\
Label Propagation                                & 0.38          & 0.36          & 0.54          & 0.75          & 0.87          & 0.90          & 0.39          & 0.65          & 0.51          \\
DORF                                             & 0.59          & 0.57          & 0.69          & \textbf{0.93} & \textbf{0.94} & 0.91          & 0.67          & \textbf{0.73} & 0.68          \\
NDORF                                            & 0.59          & 0.57          & 0.69          & \textbf{0.93} & \textbf{0.94} & 0.91          & 0.68          & \textbf{0.73} & 0.68          \\
NDSRF                                            & 0.59          & 0.57          & 0.68          & \textbf{0.93} & \textbf{0.94} & 0.91          & 0.68          & \textbf{0.73} & 0.68          \\
Rho                                              & \textbf{0.77} & 0.68          & 0.82          & 0.89          & 0.92          & 0.90          & 0.64          & 0.72          & 0.63          \\
RhoN                                             & \textbf{0.77} & 0.68          & 0.84          & 0.89          & 0.93          & \textbf{0.92} & 0.69          & 0.72          & \textbf{0.95}          \\
one\_evol                                         & 0.59          & 0.57          & \textbf{0.87} & 0.90          & 0.93          & 0.91          & \textbf{0.70} & 0.72          & \textbf{0.95} \\
qn1\_evol                                         & 0.59          & 0.57          & \textbf{0.87} & 0.90          & 0.93          & 0.91          & \textbf{0.70} & 0.72          & \textbf{0.95} \\
two\_evol                                         & 0.38          & 0.49          & 0.85          & 0.90          & 0.93          & 0.91          & \textbf{0.70} & 0.72          & \textbf{0.95} \\
qn2\_evol                                         & 0.38          & 0.49          & 0.85          & 0.90          & 0.93          & 0.91          & \textbf{0.70} & 0.72          & \textbf{0.95} \\
\botrule
\end{tabular}}
\end{table}

\begin{figure}[htbp!]
   \centering
    \includegraphics[width=1\textwidth]{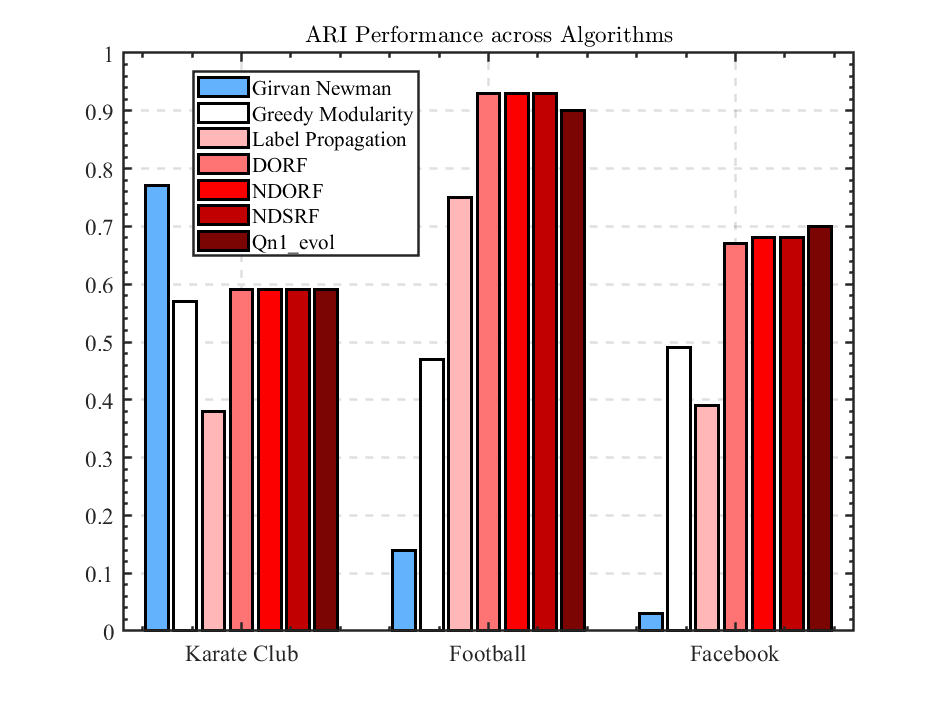}
    \caption{ARI Performance across Algorithms}
    \label{ARI Performance across Algorithms}
\end{figure}

\begin{figure}[htbp!]
   \centering
    \includegraphics[width=1\textwidth]{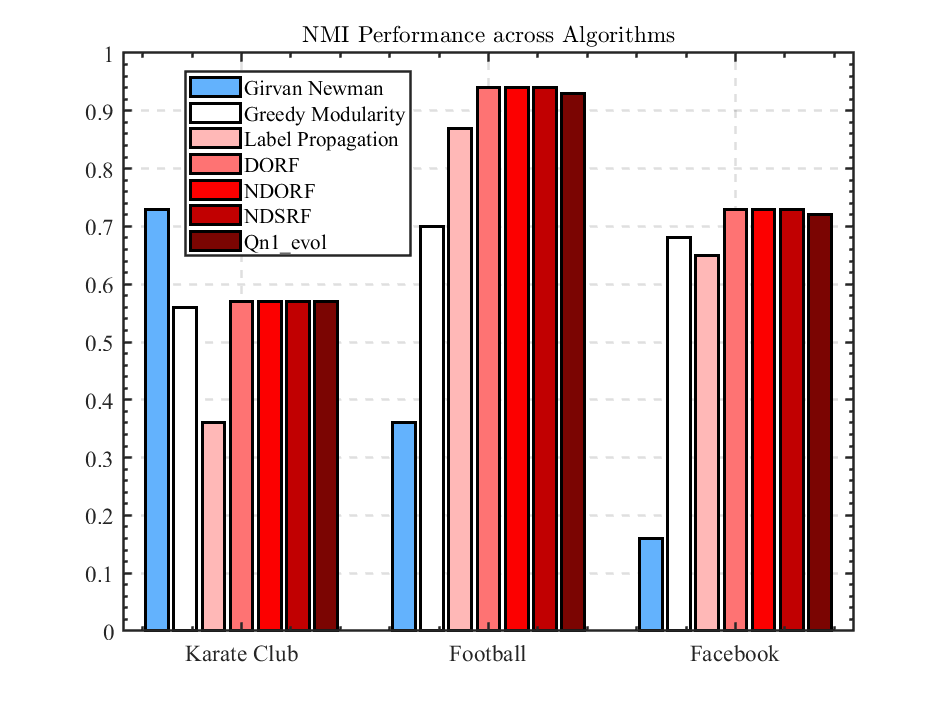}
    \caption{NMI Performance across Algorithms}
    \label{NMI Performance across Algorithms}
\end{figure}

\begin{figure}[htbp!]
   \centering
    \includegraphics[width=1\textwidth]{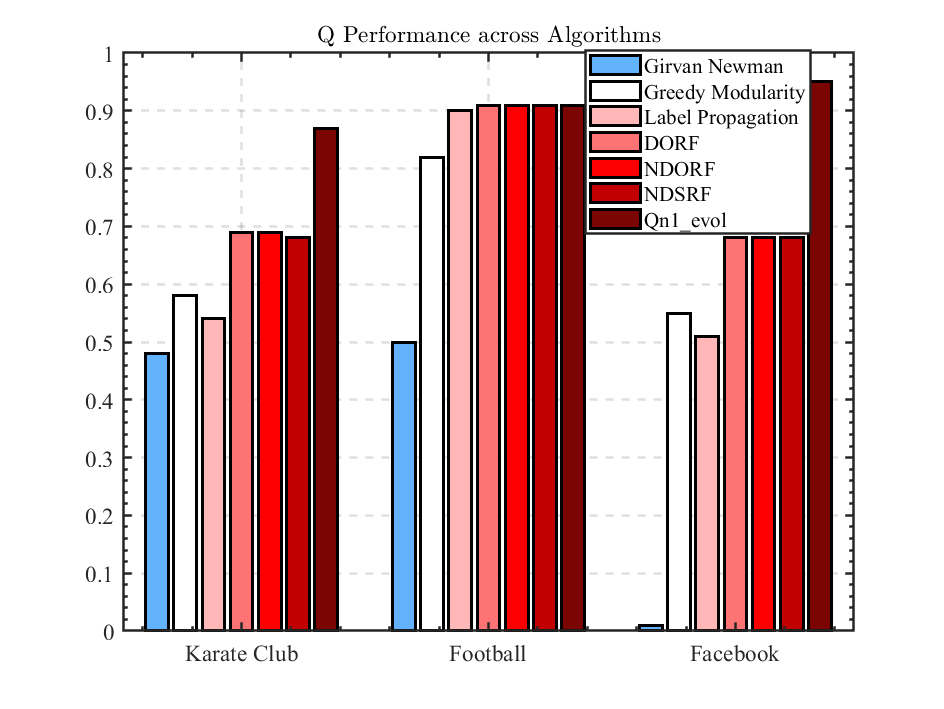}
    \caption{Q Performance across Algorithms}
    \label{Q Performance across Algorithms}
\end{figure}

The comparative performance of various algorithms on real-world networks, as depicted in the above Table \ref{tab:3} and Figures \ref{ARI Performance across Algorithms} to \ref{Q Performance across Algorithms} below, underscores the distinct advantages of evolutionary algorithms, particularly in large scale networks. Algorithms such as one\_evol, qn1\_evol, two\_evol, and qn2\_evol consistently achieve superior results in modularity, as seen in Figure \ref{Q Performance across Algorithms}, reaching a peak modularity value of 0.95 on the Facebook network.  This performance highlights evolutionary algorithms’ efficacy in optimizing community structures, especially where well-defined, connected communities are essential.

Traditional algorithms like Girvan Newman, though effective on smaller networks such as Karate Club (ARI = 0.77, NMI = 0.73), show limited scalability. For instance, on the Facebook network, its ARI drops to 0.03, and modularity Q falls to 0.01, underscoring its difficulty in managing the complexity of larger social networks. Similarly, though algorithms of DORF, NDORF, and NDSRF perform robustly on medium networks like Football (ARI = 0.93, NMI = 0.94), they fall short of the modularity optimization achieved by our evolutionary methods on larger network Facebook.

Evolutionary algorithms excel due to their ability to balance multiple performance metrics, consistently achieving the highest modularity values without sacrificing accuracy on smaller networks. On the Facebook and Football networks, they maintain strong ARI and NMI values alongside top modularity scores, indicating both accurate community detection and quality optimization of partitions. Despite their less optimal performance on smaller networks like Karate Club, evolutionary algorithms remain resilient across varying network sizes, demonstrating adaptability and superior optimization for large-scale, real-world network analysis.

 In summary, traditional methods such as Girvan Newman and those based on Ricci curvature, like the DORF series, Rho, and RhoN, have their strengths in specific contexts. However, evolutionary algorithms stand out as versatile and effective option, particularly for complex networks. Their consistent performance in achieving the highest modularity scores across various networks, along with competitive accuracy, makes them the preferred choice for network analysis when the primary goal is to optimize community structure.

\section{Concluding remarks}
In this study, we propose an evolution system for evolving weight of edge on connected finite graph according to the difference between two distances. One is the Wasserstein distance between two probability measures, the other is the distance between two vertices on the edge. Moreover, the corresponding initial value problem has been proven to have a unique global solution. Discrete version of this kind of evolution system can provide more effective algorithms for community detection than many algorithms in the same topic. Note that our systems do not require Ricci curvature on the graph. Extensive experiments on real-world datasets demonstrate that our algorithm is both easy to implement and robust across a range of parameter choices. This robustness underscores the practical applicability of our approach in accurately identifying community structures with minimal tuning. Future work may extend this method to larger, more complex networks, further enhancing its utility across diverse applications in network science.

\section*{Acknowledgements}
This research was partly supported by Public Computing Cloud, Renmin University of China.\\

\noindent
\section*{Author Contributions}
\noindent
{\bf Jicheng Ma}: Conceptualization (equal); Formal analysis (equal); Investigation (equal); Methodology (equal); Writing-original draft
(equal). {\bf Yunyan Yang}: Formal analysis (equal); Investigation (equal); Methodology (equal); Supervision (equal); Writing-review and editing
 (equal).

\section*{Data availability}
All data needed are available freely. One can find the codes of our algorithms at
https://github.com/mjc191812/Evolution-of-weights-on-a-connected-finite-graph.

\section*{Declarations}

\noindent
\textbf{Conflict of interest} The authors declared no potential conflicts of interest with respect to the research, authorship, and publication of this article.\\

\noindent
\textbf{Ethics approval} The research does not involve humans and/or animals. The authors declare that there are no ethics issues to be approved or disclosed.\\


\end{sloppypar}
\end{document}